\documentclass[twoside,12pt]{article}
\usepackage{amsmath}
\usepackage{amssymb}
\usepackage[usenames,dvipsnames]{color}
\definecolor{rosso}{rgb}{0.8,0,0}

\title{Convergence of {C}ahn--{H}illiard systems \\ 
to the {S}tefan problem\\
with dynamic boundary conditions}

\author{Takeshi Fukao\\
Department of Mathematics, Faculty of Education\\
Kyoto University of Education\\
1~Fujinomori, Fukakusa, Fushimi-ku, Kyoto~612-8522 Japan\\
E-mail: \texttt{fukao@kyokyo-u.ac.jp}}
 
\date{}

\newcommand\testopari{\sc Takeshi Fukao}
\newcommand\testodispari{\sc {S}tefan problem with dynamic boundary conditions}
\begin{document}

\maketitle

\begin{abstract}
This paper examines the well-posedness of the {S}tefan problem with 
a dynamic boundary condition. 
To show the existence of the weak solution, 
the original problem is approximated by a limit 
of an equation and dynamic boundary condition of {C}ahn--{H}illiard type. 
By using this {C}ahn--{H}illiard approach, 
it becomes clear that 
the state of the mushy region of the {S}tefan problem is 
characterized by an asymptotic limit of the fourth-order 
system, which has a double-well structure. 
This fact also raises the possibility of the numerical application  
of the {C}ahn--{H}illiard system to the degenerate parabolic equation, 
of which the {S}tefan problem is one. 

\vspace{2mm}

\vspace{2mm}
\noindent \textbf{AMS (MOS) subject clas\-si\-fi\-ca\-tion:} 80A22, 35K61, 35K25, 35D30, 47J35.

\end{abstract}

%%%%% Section 1. %%%%%
\section{Introduction}
\setcounter{equation}{0}

The {S}tefan problem is a well-known mathematical model 
that describes 
the solid--liquid phase transition. 
Among many results in the literature (for example, \cite{Fri68, Kam61, MVV83, Vis96} and so on),
the following enthalpy formulation of the {S}tefan problem 
with the {D}irichlet--{R}obin type boundary condition was studied by the $L^2$-framework in \cite{Dam77}: 
\begin{equation} 
	\frac{\partial u}{\partial t} -\Delta \beta (u) = g 
	\quad \mbox{in }Q:=(0,T) \times \Omega, 
	\label{ST1}
\end{equation}
where $0<T<+\infty$ is a finite time and 
$\Omega \subset \mathbb{R}^{d}$, $d=2$ or $3$, 
is a bounded domain with smooth boundary $\Gamma $.
The unknown $u:\overline{Q} \to \mathbb{R}$ denotes the enthalpy and 
$\beta (u)$ is the temperature; 
$g:Q \to \mathbb{R}$ 
is a given heat source. 
In the model of the {S}tefan problem 
$\beta :\mathbb{R} \to \mathbb{R}$ is a piecewise linear function 
of the following form: 
\begin{equation} 
	\beta (r):=
	\begin{cases}
	k_s r & r<0, \\
	0     & 0\le r \le L, \\
	k_\ell (r-L) & r >L;
	\end{cases} 
	\label{beta}
\end{equation}
$k_s$, $k_\ell >0$ represet the heat conductivities on the 
solid and liquid regions, respectively, and $L>0$ is the latent heat constant. 
Consider the initial-boundary value problem for this kind of 
partial differential equation. 
The dynamic boundary condition is a sort of differential equation that 
represents the dynamics on the boundary $\Gamma$. 
As the condition for solving the partial differential equation, the 
time derivative is included, and this is well treated like, {D}irichlet, {N}eumann, and 
{R}obin-type conditions for various problems. 
Under the dynamic boundary condition of the following form:
\begin{equation} 
	\frac{\partial \beta (u)}{\partial t} + \partial_{\mbox{\scriptsize \boldmath $ \nu $}} \beta (u)  = g_\Gamma 
	\quad \mbox{on }\Sigma:=(0,T) \times \Gamma,
	\label{aiki}
\end{equation}
the existence and uniqueness of \eqref{ST1} was 
studied in \cite{Aik95, Aik96} using a subdifferential approach, 
where the symbol  
$\partial_{\mbox{\scriptsize \boldmath $ \nu $}} $ 
denotes the normal derivative on $\Gamma $ outward from $\Omega$;
$g_\Gamma: \Sigma \to \mathbb{R}$ is a given heat source on the boundary. 
For a more general setting, we can find the result in \cite{AMTI06}.

In this paper, 
the well-posedness of the {S}tefan problem with the following
dynamic boundary condition is studied:
\begin{equation} 
	\frac{\partial u}{\partial t} + \partial_{\mbox{\scriptsize \boldmath $ \nu $}} \beta (u) -\Delta _\Gamma \beta (u) = g_\Gamma 
	\quad \mbox{on }\Sigma,
	\label{fukao}
\end{equation}
where the symbol $\Delta _{\Gamma }$ stands for the {L}aplace--{B}eltrami operator
on $\Gamma $ (see, e.g., \cite[Chapter~3]{Gri09}). 
If we simultaneously consider \eqref{ST1} 
on the bulk $\Omega$ and \eqref{aiki} or \eqref{fukao} 
on the boundary $\Gamma$, 
then the setting of \eqref{fukao} is more natural than that of \eqref{aiki}. 

The main idea of the existence result is 
to approximate the {S}tefan problem from 
the {C}ahn--{H}illiard system. 
In 2009, {G}oldstein, {M}iranville, and {S}chimperna studied the following  
{C}ahn--{H}illiard system: For $\varepsilon > 0$
\begin{gather} 
	\frac{\partial u}{\partial t} -\Delta \mu = 0 
	\quad \mbox{in }Q, 
	\label{CH1}
	\\
	\mu = -\varepsilon \Delta u + \beta (u)+ \varepsilon \pi (u)-f 
	\quad \mbox{in }Q,
	\label{CH2}
\end{gather}
with a dynamic boundary condition of the following form:
\begin{gather} 
	\frac{\partial u}{\partial t} + \partial_{\mbox{\scriptsize \boldmath $ \nu $}} \mu - \Delta _\Gamma \mu  =0
	\quad \mbox{on }\Sigma,
	\label{CHB1}
	\\
	\mu =\varepsilon\partial _{\mbox{\scriptsize \boldmath $ \nu $}} u - \varepsilon\Delta _\Gamma u 
	+ \beta _\Gamma (u)
	+ \varepsilon\pi _\Gamma (u)-f_\Gamma 
	\quad \mbox{on }\Sigma.
	\label{CHB2}
\end{gather}
The unknowns $u, \mu :\overline{Q} \to \mathbb{R}$ represent the 
order parameter and chemical potential, respectively. 
Let us recall some basic concepts.
The {C}ahn--{H}illiard system is characterized by 
the nonlinear terms $\beta +\pi $ and $\beta _\Gamma +\pi _\Gamma $, 
which are some derivatives 
of functions $W$ and $W_\Gamma $, respectively. 
Usually referred to as double-well potentials, 
for example, $W(r)=W_\Gamma (r)=(1/4)(r^2-1)^2$, in this 
case $\beta (r)=\beta _\Gamma (r)=r^3$ and 
$\pi(r)=\pi _\Gamma (r)=-r$ for all $r \in \mathbb{R}$. 
For the details, see \cite{CH58, EZ86}. 
Now, taking $\beta _\Gamma =\beta $ as \eqref{beta} and 
letting $\varepsilon \to 0$ in \eqref{CH2}, we obtain \eqref{ST1} 
as the limiting problem of \eqref{CH1} with $\mu =\beta(u)-f$ 
and some 
modification of $-\Delta f=g$ in $\Omega$ and 
$\partial _{\mbox{\scriptsize \boldmath $ \nu $}} f-\Delta_\Gamma f=g_\Gamma $ in 
$\Gamma$, with the setting that the trace $f_{|_\Gamma }$ of $f$ coincides with $f_\Gamma$. 
Then, from \eqref{CHB2}, we can also characterize 
the dynamic boundary condition \eqref{fukao} by \eqref{CHB1}.

The existence and uniqueness of 
the {C}ahn--{H}illiard system \eqref{CH1}--\eqref{CHB2}, as well as its asymptotic behavior,
was treated by \cite{GMS11} 
in the case when $\varepsilon >0$. 
Following this, many related problems were treated in \cite{CGM13, CP14}. 
Recently, the author extended the pioneering work of \cite{GMS11} to 
the more general case in which
$\beta $ and $\beta _\Gamma $ are maximal monotone 
graphs. This included the singular case for the 
subdifferential $\partial I_{[-1,1]}$ of the 
indicator function $I_{[-1,1]}$ of the closed interval $[-1,1]$, 
and 
when $\pi $ and $\pi _\Gamma $ are Lipschitz continuous functions. 
The well-posedness of strong and weak solutions for the 
initial-boundary problem 
under appropriate assumptions was discussed in \cite{CF15}, where 
the essential idea comes from 
\cite{CC13, CGS14}.

The present paper proceeds as follows.
In Section~2, the main theorem is stated. 
First, we prepare the notation used in this paper.
Next, we recall a known result for some equations and dynamic boundary conditions of 
{C}ahn--{H}illiard type (P)$_\varepsilon$ in Proposition~2.1. 
The main theorem is related to the
convergence of {C}ahn--{H}illiard systems (P)$_\varepsilon$ to the {S}tefan problem (P)
with dynamic boundary conditions. 
In Section~3, we obtain a uniform estimate that is useful in proving the main theorem. 
To guarantee sufficient regularity of the unknowns for (P)$_\varepsilon$, 
we start from the approximate problem (P)$_{\varepsilon, \lambda}$. 
After obtaining all necessary estimates, we correct similar uniform estimates 
for (P)$_\varepsilon$. 
A proof of the main theorem is given in Section~4. The 
strategy of the proof proceeds in a standard manner. Based on the uniform estimates, 
we consider the limiting procedure $\varepsilon \to 0$. 
The main theorem is applied under a more general assumption for $\beta $, namely that it is
not only a non-decreasing piecewise linear function, as in \eqref{beta}, but also 
some maximal monotone graph. 
Therefore, we apply the monotonicity argument. 
Finally, in Section~5, the uniqueness of (P) is proved by showing 
the continuous dependence for the given data. 

A detailed index of sections and subsections is as follows:

\begin{itemize}
 \item[1.] Introduction
 \item[2.] Main results
	\begin{itemize}
	 \item[2.1.] Notation
	 \item[2.2.] Main theorem
	 \item[2.3.] Approximate solutions to the {C}ahn--{H}illiard system (P)$_\varepsilon$
	\end{itemize}
 \item[3.] Uniform estimates
	\begin{itemize}
	 \item[3.1] Uniform estimates for approximate solutions of (P)$_{\varepsilon,\lambda}$
	 \item[3.2] Uniform estimates for approximate solutions of (P)$_\varepsilon$
	\end{itemize}
 \item[4.] Proof of the main theorem
 \item[5.] Continuous dependence
\end{itemize}

%%%%% Section 2. %%%%%
\section{Main results}
\setcounter{equation}{0}

In this section, the main theorem is stated. 
First, 
we recall a previous 
result from \cite{CF15} that plays an important role in this paper. 
Next, we present the main theorem.

%%%%% Section 2.1. %%%%%
\subsection{Notation}
We use the 
spaces $H:=L^2(\Omega )$, $V:=H^1(\Omega )$, 
$H_\Gamma :=L^2(\Gamma )$, and $V_\Gamma :=H^1(\Gamma )$
with the usual norms 
$| \cdot |_{H}$, $|\cdot |_{V}$, $|\cdot |_{H_\Gamma}$, $|\cdot |_{V_\Gamma}$ 
and inner products $(\cdot,\cdot )_{H}$, $(\cdot ,\cdot )_{V}$,
$(\cdot,\cdot )_{H_\Gamma}$, $(\cdot ,\cdot )_{V_\Gamma}$, respectively. 
Moreover, $\mbox{\boldmath $ H$}:=H \times H_\Gamma$, 
$\mbox{\boldmath $ V$}:=\{ (z,z_\Gamma ) \in V \times V_\Gamma 
: z_\Gamma =z_{|_\Gamma} \mbox{ a.e.\ on }\Gamma  \}$, and 
$\mbox{\boldmath $ W$}:=H^2(\Omega ) \times H^2(\Gamma )$. Hereafter, we use a bold symbol $\mbox{\boldmath $ z$}$ to denote the pair 
$(z,z_\Gamma )$ corresponding to the letter. 
Then, $\mbox{\boldmath $ H$}$, $\mbox{\boldmath $ V$}$, and $\mbox{\boldmath $ W$}$ are 
{H}ilbert spaces with 
the inner product
\begin{equation*}
	(\mbox{\boldmath $ u $},\mbox{\boldmath $ z $}
	)_{\mbox{\scriptsize \boldmath $ H$}}
	:=(u,z)_{H} + (u_\Gamma ,z_\Gamma )_{H_\Gamma } \quad 
	\mbox{for all}~\mbox{\boldmath $ u$}, 
	\mbox{\boldmath $ z$}
	\in \mbox{\boldmath $ H$},
\end{equation*}
and the related norm 
is analogously defined as one of $\mbox{\boldmath $ V$}$ or $\mbox{\boldmath $ W$}$. Note that, if 
$\mbox{\boldmath $ z$} \in \mbox{\boldmath $ V$}$, then $z_{\Gamma }$ 
is exactly the trace $z_{|_\Gamma}$ of $z$ on $\Gamma$, 
whereas if $\mbox{\boldmath $ z$}$ is just in $ \mbox{\boldmath $ H$}$, 
then $z \in H$ and $z_{\Gamma } \in H_{\Gamma }$ are independent. 
Define $m:\mbox{\boldmath $ H$} \to \mathbb{R}$ by 
\begin{equation}
	m(\mbox{\boldmath $ z$}):=\frac{1}{|\Omega |+|\Gamma| }
	\left\{ \displaystyle \int_{\Omega }^{}z dx
	+ \int_{\Gamma }^{} z_{\Gamma } d\Gamma \right\} 
	\quad \mbox{for all }\mbox{\boldmath $ z$} \in \mbox{\boldmath $ H$},
\label{mean}
\end{equation}
where $|\Omega |:=\int_{\Omega }^{}1dx$ and $|\Gamma |:=\int_{\Gamma }^{}1d\Gamma $. 
The symbol $\mbox{\boldmath $ V$}^*$ denotes 
the dual space of $\mbox{\boldmath $ V$}$, 
and the pair 
$\langle \cdot ,\cdot 
\rangle _{\mbox{\scriptsize \boldmath $ V$}^*, \mbox{\scriptsize 
\boldmath $ V$}}$
denotes the duality pairing between $\mbox{\boldmath $ V$}^*$ and 
$\mbox{\boldmath $ V$}$. 
Moreover, we define 
the bilinear form 
$a(\cdot ,\cdot ):\mbox{\boldmath $ V$} \times \mbox{\boldmath $ V$} \to \mathbb{R}$ by 
\begin{equation*}
	a(\mbox{\boldmath $ u$},\mbox{\boldmath $ z$} ):=
	\int_{\Omega }^{} \nabla u \cdot \nabla z dx 
	+\int_{\Gamma }^{} \nabla _\Gamma u _\Gamma \cdot \nabla _\Gamma z_\Gamma d\Gamma 
	\quad \mbox{for all }\mbox{\boldmath $ u$},\mbox{\boldmath $ z$} \in \mbox{\boldmath $ V$},
\end{equation*} 
where 
$\nabla _{\Gamma }$ denotes the surface gradient on $\Gamma$ (see, e.g., \cite[Chapter 3]{Gri09}). 
We also introduce 
the subspace of $\mbox{\boldmath $ H$}$ as 
$\mbox{\boldmath $ H$}_0:=\{ \mbox{\boldmath $ z$} \in
\mbox{\boldmath $ H$} : m(\mbox{\boldmath $ z$})=0 \}$
and $\mbox{\boldmath $ V$}_0 :=\mbox{\boldmath $ V$} \cap \mbox{\boldmath $ H$}_0$ with 
their norms $| \mbox{\boldmath $ z$}|_{\mbox{\scriptsize \boldmath $ H$}_0}:=|\mbox{\boldmath $ z$}|_{\mbox{\scriptsize \boldmath $ H$}}$ for all 
$\mbox{\boldmath $ z$} \in \mbox{\boldmath $ H$}_0$ and 
$|\mbox{\boldmath $ z$}|_{\mbox{\scriptsize \boldmath $ V$}_0}:=
a(\mbox{\boldmath $ z$},\mbox{\boldmath $ z$} )^{1/2}
$
for all 
$\mbox{\boldmath $ z$} \in \mbox{\boldmath $ V$}_0$. 
Then, we can define the duality mapping 
$\mbox{\boldmath $ F$}: \mbox{\boldmath $ V$}_0 \to \mbox{\boldmath $ V$}_0^*$ by 
\begin{equation}
	\langle \mbox{\boldmath $ F$} 
	\mbox{\boldmath $ z$}, \tilde{\mbox{\boldmath $ z$}} 
	\rangle _{\mbox{\scriptsize \boldmath $ V$}_0^*, 
	\mbox{\scriptsize \boldmath $ V$}_0}
	:= a(\mbox{\boldmath $ z$},\tilde{\mbox{\boldmath $ z$}})
	\quad 
	{\rm for~all~}
	\mbox{\boldmath $ z$}, \tilde{\mbox{\boldmath $ z$}} \in \mbox{\boldmath $ V$}_0.
\label{dual}
\end{equation}
Using this, we can define the inner product in $\mbox{\boldmath $ V$}_0^*$ by 
\begin{equation}
	(\mbox{\boldmath $ z$}_1^*,\mbox{\boldmath $ z$}_2^*)_{\mbox{\scriptsize \boldmath $ V$}_0^*}
	:=\langle \mbox{\boldmath $ z$}_1^*, 
	\mbox{\boldmath $ F$} ^{-1} \mbox{\boldmath $ z$}_2^* 
	\rangle _{\mbox{\scriptsize \boldmath $ V$}^*_0, \mbox{\scriptsize \boldmath $ V$}_0}
	\quad \mbox{for all } \mbox{\boldmath $ z$}_1^*,\mbox{\boldmath $ z$}_2^* \in \mbox{\boldmath $ V$}_0^*.
	\label{Vstar}
\end{equation}
Then, we obtain the dense and compact embeddings
$\mbox{\boldmath $ V$}_0 
\mathop{\hookrightarrow} \mathop{\hookrightarrow}
\mbox{\boldmath $ H$}_0 
\mathop{\hookrightarrow} \mathop{\hookrightarrow}
\mbox{\boldmath $ V$}_0^*$; see \cite{CF15} for details of this setting. 
These are essentially the same as in previous studies \cite{CF14, KN96, Kub12}.

%%%%% Section 2.3. %%%%%
\subsection{Main theorem}

In this subsection, 
we define the weak solution for the {S}tefan problem with a dynamic boundary condition. 
Then, we give the main theorem.

First, we present the target problem (P) 
of this paper, namely the {S}tefan problem with 
the dynamic boundary condition: 
\begin{gather} 
	\frac{\partial u}{\partial t}-\Delta \beta (u) = g 
	\quad \mbox{\rm a.e.\ in }Q,
	\nonumber 
%	\label{Sae1}
\\
	\beta (u_\Gamma) = \beta (u)_{|_\Gamma }, 
	\quad 
	\frac{\partial u_\Gamma }{\partial t}
	+ \partial_{\mbox{\scriptsize \boldmath $ \nu $}} \beta (u)
	- \Delta _\Gamma \beta (u_\Gamma) =g_\Gamma 
	\quad \mbox{\rm a.e.\ on }\Sigma,
	\label{Sae2}
\\
	u(0)=u_0
	\quad 
	\mbox{\rm a.e.\ in }\Omega, \quad 
	u_\Gamma  (0)=u_{0\Gamma}
	\quad \mbox{\rm a.e.\ on }\Gamma,
	\nonumber 
%	\label{Sae3}
\end{gather} 
where the prototype {S}tefan problem 
is formulated by the setting \eqref{beta}. 
In this paper, the target problem (P) will be 
formulated in a more general setting; see Definition~2.1.

\paragraph{Remark 1.} The dynamic boundary condition 
is arranged from three previous results 
\cite{Aik95, Aik96, AMTI06} regarding not only 
the {L}aplace--{B}eltrami operator, 
but also the time derivative. 
Actually, they treated $\partial \beta (u_\Gamma) /\partial t$. 
In the case of \eqref{beta}, 
the first condition \eqref{Sae2} implies that 
$u_\Gamma $ is not necessarily equal to the trace $u_{|_\Gamma }$ of $u$. More 
precisely, we will obtain 
$u_\Gamma =u_{|_\Gamma }$ except in the mushy region 
$\{ x \in \Omega : 0 \le u(x) \le L \}$
of the {S}tefan problem. 
This is because we can only expect $\beta (u(t)) \in V$, 
but $u(t) \not\in V$ (weak solution). \\

Throughout this paper, we assume that 
\begin{itemize}
 \item[(A1)] $\beta$ is a maximal monotone graph in $\mathbb{R}\times \mathbb{R}$, 
and is a subdifferential $\beta =\partial \widehat{\beta }$ 
of some proper, lower semicontinuous, and convex function 
$\widehat{\beta}:\mathbb{R} \to [0,+\infty ]$
satisfying $\widehat{\beta }(0)=0$ with some effective domain $D(\beta )$. This implies 
$\beta (0)=0$. Moreover, 
there exist two constants $c_1$, $c_2>0$ such that 
\begin{equation} 
	\widehat{\beta} (r) \ge c_1 |r|^2- c_2 \quad {\rm for~all}~r \in \mathbb{R};
	\label{prim0}
\end{equation}
 \item[(A2)] $\pi :D(\pi)=\mathbb{R} \to \mathbb{R}$ is 
{L}ipschitz continuous. Moreover, 
\begin{equation} 
	|\pi '(r)| \le 1 \quad {\rm for~a.a.~} r \in \mathbb{R};
	\label{pi}
\end{equation} 
 \item[(A3)] $\mbox{\boldmath $ g$} \in L^2(0,T;\mbox{\boldmath $ H$}_0)$; 
 \item[(A4)] $\mbox{\boldmath $ u$}_0:=(u_0,u_{0\Gamma}) \in \mbox{\boldmath $ V$}$ with $m_0 \in {\rm int}D(\beta )$, and 
the compatibility conditions $\widehat{\beta }(u_0)\in L^1(\Omega ), 
\widehat{\beta }(u_{0\Gamma }) \in L^1(\Gamma)$ hold. 
\end{itemize}

For simplicity, we have assumed that the derivative of $\pi $ is bounded by $1$ in 
\eqref{pi}. It is sufficient to assume that $\pi \in W^{1,\infty }(\mathbb{R})$. 
We now define a weak solution.

\paragraph{Definition 2.1.} 
{\it The pair $(\mbox{\boldmath $ u$}, \mbox{\boldmath $ \xi $})$ of 
functions $\mbox{\boldmath $ u$} \in H^1(0,T;\mbox{\boldmath $ V$}^*) \cap L^2(0,T;\mbox{\boldmath $ H$})$ and 
$\mbox{\boldmath $ \xi $} \in L^2(0,T;\mbox{\boldmath $ V$})$
is called the weak solution of {\rm (P)} if} 
\begin{equation*} 
	\xi \in \beta (u) \quad \mbox{a.e.\ in } Q, 
	\quad 
	\xi _\Gamma \in \beta (u_\Gamma ), \quad 
	\xi _\Gamma =\xi _{|_\Gamma } \quad \mbox{a.e.\ on } \Sigma, 
\end{equation*} 
{\it and they satisfy}
\begin{align} 
	\bigl \langle u' (t),z 
	\bigr \rangle _{V^*,V}
	& + \langle u'_{\Gamma} (t),z_\Gamma 
	\bigr \rangle _{V_\Gamma ^*,V_\Gamma }
	+ \int_{\Omega }^{} \nabla \xi (t) \cdot \nabla z dx 
	+ \int_{\Gamma }^{} \nabla_\Gamma  \xi _{\Gamma}(t) 
	\cdot \nabla _\Gamma z_\Gamma 
	d \Gamma 
	\nonumber \\
	& 
	= \int_{\Omega }^{} g(t) z dx 
	+ \int_{\Gamma }^{} g_\Gamma (t) z_\Gamma 
	d \Gamma 
	\quad {\it for~all~}\mbox{\boldmath $ z$}=(z,z_\Gamma ) \in \mbox{\boldmath $ V$}
	\label{def}
\end{align}
{\it for a.a.\ $t \in (0,T)$ with 
$u(0)=u_0$
a.e.\ in $\Omega$ and 
$u_\Gamma  (0)=u_{0\Gamma}$
a.e.\ on $\Gamma$.} \\

Our main theorem is now stated.

\paragraph{Theorem 2.1.} {\it Assume {\rm (A1)}--{\rm (A4)} hold. Then, 
there exists at least one weak solution 
$( \mbox{\boldmath $ u$}, \mbox{\boldmath $ \xi $})$ of {\rm (P)}.} \\

Moreover, we obtain the 
following continuous dependence for the given data:

\paragraph{Theorem 2.2.} 
{\it Assume {\rm (A1)}--{\rm (A4)}. For $i=1,2$, 
let $(\mbox{\boldmath $ u$}^{(i)}, 
\mbox{\boldmath $ \xi $}^{(i)})$ be a weak solution of {\rm (P)}  
corresponding to the data $\mbox{\boldmath $ g$}^{(i)}$ and 
$\mbox{\boldmath $ u$}_0^{(i)}$, in which we assume that 
$m(\mbox{\boldmath $ u$}_0^{(1)})=m(\mbox{\boldmath $ u$}_0^{(2)})$. 
Then, there exists a positive constant $C$ that depends only on $T$ such that
\begin{equation} 
	\bigl| \mbox{\boldmath $ u$}^{(1)}(t)-  \mbox{\boldmath $ u$}^{(2)}(t) 
	\bigr|_{\mbox{\scriptsize \boldmath $ V$}_0^*}^2 
	\le C \left\{ 
	\bigl| \mbox{\boldmath $ u$}^{(1)}_0-  \mbox{\boldmath $ u$}^{(2)}_0 
	\bigr|_{\mbox{\scriptsize \boldmath $ V$}_0^*}^2 
	+
	\int_{0}^{T} 
	\bigl| \mbox{\boldmath $ g$}^{(1)}(s ) -\mbox{\boldmath $ g$}^{(2)}(s ) 
	\bigr|_{\mbox{\scriptsize \boldmath $ H$}_0}^2 ds
	\right\}
	\label{conti1}
\end{equation} 
for all $t\in [0,T]$. 
Moreover, if $\beta $ is Lipschitz continuous, then
\begin{align} 
	\lefteqn{ 
	\int_{0}^{t} 
	\bigl| \mbox{\boldmath $ \xi$}^{(1)}(s)-  \mbox{\boldmath $ \xi $}^{(2)}(s) 
	\bigr|_{\mbox{\scriptsize \boldmath $ H$}}^2 
	ds
	} \nonumber \\
	& \le \frac{c_\beta }{2} C( 1+T)  \left\{ 
	\bigl| \mbox{\boldmath $ u$}^{(1)}_0-  \mbox{\boldmath $ u$}^{(2)}_0 
	\bigr|_{\mbox{\scriptsize \boldmath $ V$}_0^*}^2 
	+
	\int_{0}^{T} 
	\bigl| \mbox{\boldmath $ g$}^{(1)}(s ) -\mbox{\boldmath $ g$}^{(2)}(s ) 
	\bigr|_{\mbox{\scriptsize \boldmath $ H$}_0}^2 ds
	\right\}
	\label{conti2}
\end{align}
for all $t \in [0,T]$, 
where $c_\beta >0$ is the {L}ipschitz constant of $\beta $. }\\

Note that 
the assumption 
$m(\mbox{\boldmath $ u$}_0^{(1)})=m(\mbox{\boldmath $ u$}_0^{(2)})$ gives us 
that 
$\mbox{\boldmath $ u$}^{(1)}(t)-\mbox{\boldmath $ u$}^{(2)}(t)
\in \mbox{\boldmath $ V$}_0^*$ 
for all $t \in [0,T]$.
This theorem implies the uniqueness of the weak solution of (P).

%%%%% Section 2.3. %%%%%
\subsection{Approximate solutions to the {C}ahn--{H}illiard system (P)$_\varepsilon$}

In this subsection, 
we state the approximate problem for (P). 
For this aim, we recall a previous result \cite{CF15} for 
the equation and dynamic boundary condition 
of the 
{C}ahn--{H}illiard type
(P)$_\varepsilon $.  
This can be written as the following initial-boundary value problem \eqref{CHae1}--\eqref{CHae5}:
For $\varepsilon >0$,
\begin{gather} 
	\frac{\partial u}{\partial t}-\Delta \mu = 0 
	\quad \mbox{\rm a.e.\ in }Q,
	\label{CHae1}
\\
	\xi \in \beta (u), \quad \mu = -\varepsilon \Delta u + \xi + \varepsilon \pi (u) -f 
	\quad \mbox{a.e.\ in }Q,
	\label{CHae2}
\\
	u_\Gamma =u _{|_\Gamma }, \quad \mu _\Gamma =\mu _{|_\Gamma },
	\quad 
	\frac{\partial u_\Gamma }{\partial t}
	+ \partial_{\mbox{\scriptsize \boldmath $ \nu $}} \mu 
	- \Delta _\Gamma \mu _\Gamma =0
	\quad \mbox{a.e.\ on }\Sigma,
	\label{CHae3}
\\ 
	\xi _\Gamma \in \beta (u_\Gamma ), \quad \mu _\Gamma =\varepsilon \partial_{\mbox{\scriptsize \boldmath $ \nu $}} u - 
	\varepsilon \Delta _\Gamma u_\Gamma +\xi _\Gamma + \varepsilon \pi  (u_\Gamma )-f_\Gamma 
	\quad \mbox{a.e.\ on }\Sigma,
	\label{CHae4}
\\
	u(0)=u_0
	\quad 
	\mbox{\rm a.e.\ in }\Omega, \quad 
	u_\Gamma  (0)=u_{0\Gamma}
	\quad \mbox{\rm a.e.\ on }\Gamma,
	\label{CHae5}
\end{gather} 
where 
$f : Q \to \mathbb{R}$, $f_\Gamma : \Sigma \to \mathbb{R}$, 
$u_0 : \Omega \to \mathbb{R}$, and $u_{0\Gamma } : \Gamma \to \mathbb{R}$ are 
given data. In particular, $f$ and $f_\Gamma$ are constructed by $g$ and $g_\Gamma$ at a later point. 
In the main theorem, we treated the general setting (A1) and (A2) of 
$\beta $ for the 
{S}tefan problem with some suitable $\pi$. 
In the setting of \eqref{beta},
one example of  
$\pi :\mathbb{R} \to \mathbb{R}$ is a piecewise linear function 
of the following form: 
\begin{equation} 
	\beta (r):=
	\begin{cases}
	k_s r & {\rm if}~r<0, \\
	0     & {\rm if}~0\le r \le L, \\
	k_\ell (r-L) & {\rm if}~r >L,
	\end{cases} 
	\quad 
	\pi (r):=
	\begin{cases}
	\displaystyle \frac{L}{2} & {\rm if}~r<0, \vspace{2mm}\\
	\displaystyle \frac{L}{2}-r     & {\rm if}~0\le r \le L, \vspace{2mm}\\
	\displaystyle -\frac{L}{2} & {\rm if}~r >L,
	\end{cases} 
	\label{betapi}
\end{equation}
of course, in this case $\xi =\beta (u)$ and $\xi _\Gamma =\beta (u_\Gamma)$.  
Therefore, (A1) and (A2) hold; 
actually $\widehat{\beta }$ is obtained by 
\begin{equation*} 
	\widehat{\beta} (r):=
	\begin{cases}
	\displaystyle \frac{k_s}{2} r^2 & {\rm if}~r<0, \vspace{2mm}\\
	0     & {\rm if}~0\le r \le L, \vspace{2mm}\\
	\displaystyle \frac{k_\ell}{2} (r-L)^2 & {\rm if}~r >L,
	\end{cases} 
	\label{betahat}
\end{equation*}
so that $\beta =\partial \widehat{\beta}$, and this quadratic function $\widehat{\beta }$ satisfies \eqref{prim0}.
Because of this value of $\varepsilon \pi $,
we can realize the double-well structure of 
the potential 
$W=\widehat{\beta }+\varepsilon \widehat{\pi}$ as the sum of primitives of 
$\beta $ and $\varepsilon \pi $. Therefore, 
(P)$_\varepsilon $ with the prototype setting 
\eqref{betapi} has the exact structure 
of the {C}ahn--Hilliard system.

\paragraph{Remark 2.} 
Viewed in terms of the {S}tefan problem, $\pi $ is an artificial term. 
Therefore, we can prove the existence of the weak 
solution of the {S}tefan problem even if $\pi \equiv 0$. 
However, in this case, we no longer find the 
structure of the {C}ahn--{H}illiard system (more precisely, 
the double-well structure) at the 
level of (P)$_\varepsilon$. 
If we construct the relationship between 
some equation of {C}ahn--{H}illiard type and the original {S}tefan problem, 
we must choose a suitable $\pi $ depending on $\beta $, 
which breaks the monotonicity of $\beta $. 
This treatment is completely independent of the choice of boundary condition. 
We will focus on the convergence from the {C}ahn--{H}illiard system 
to the nonlinear diffusion equation under the {N}eumann boundary condition in a 
forthcoming paper \cite{CF16}. 
\\

Here, we know that for each 
$\mbox{\boldmath $ g$}:=(g,g_\Gamma ) \in L^2(0,T;\mbox{\boldmath $ H$}_0)$, 
there exists $\mbox{\boldmath $ f$} \in L^2(0,T;\mbox{\boldmath $ V$}_0)$ such that 
\begin{align}
	\int_{\Omega }^{} \nabla f \cdot \nabla z dx + \int_{\Gamma }^{} 
	\nabla _{\Gamma } f_\Gamma \cdot \nabla _\Gamma z_\Gamma d\Gamma
	=  \int_{\Omega }^{} g  z dx + \int_{\Gamma }^{} 
	 g_\Gamma  z_\Gamma d\Gamma
	\quad {\rm for~all}~\mbox{\boldmath $ z$} \in \mbox{\boldmath $ V$}_0,
	\label{ax}
\end{align}
a.e.\ in $(0,T)$ (see \eqref{phi}, \eqref{star}, also \cite[Lemma~C]{CF15}). 
To state the previous known results, we 
use the change of variable $\mbox{\boldmath $ v$}
:=\mbox{\boldmath $ u$}-m_0 \mbox{\boldmath $ 1$}$, 
where 
$m_0:=m(\mbox{\boldmath $ u$}_0)$ and 
$\mbox{\boldmath $ 1$}:=(1,1)$, the pair of constants. 
We also use
$\mbox{\boldmath $ v$}_0:=\mbox{\boldmath $ u$}_0-m_0 \mbox{\boldmath $ 1$}$. 
We set $\mbox{\boldmath $ \beta $}(\mbox{\boldmath $ z$}):=(\beta (z),\beta (z_\Gamma ))$, 
$\mbox{\boldmath $ \pi $}(\mbox{\boldmath $ z$}):=(\pi (z),\pi (z_\Gamma ))$ for 
all $\mbox{\boldmath $ z$} \in \mbox{\boldmath $ H$}$. 
Under assumptions (A1)--(A4), 
\cite[Theorem 2.2]{CF15} leads to the following proposition of the 
existence and uniqueness for the equation and 
dynamic boundary condition of {C}ahn--{H}illiard type 
(P)$_\varepsilon $ with respect to our $\beta +\varepsilon \pi $.

\paragraph{Proposition 2.1.} {\it Under assumptions {\rm (A1)}--{\rm (A4)}, 
there exists a triplet $(\mbox{\boldmath $ v$}_\varepsilon, 
\mbox{\boldmath $ \mu  $}_\varepsilon, \mbox{\boldmath $ \xi $}_\varepsilon )$ 
of} 
\begin{gather*}
	\mbox{\boldmath $ v$}_\varepsilon 
	\in H^1(0,T;\mbox{\boldmath $ V$}_0^*) \cap L^\infty (0,T;\mbox{\boldmath $ V$}_0) \cap L^2(0,T;\mbox{\boldmath $ W$}), \\
	\mbox{\boldmath $ \mu $}_\varepsilon 
	\in L^2(0,T;\mbox{\boldmath $ V$}), \quad 
	\mbox{\boldmath $ \xi $}_\varepsilon 
	\in L^2(0,T;\mbox{\boldmath $ H$}),
\end{gather*}
{\it that satisfy}
\begin{gather} 
	\bigl \langle \mbox{\boldmath $ v$}'_\varepsilon (t),\mbox{\boldmath $ z$} 
	\bigr \rangle _{\mbox{\scriptsize \boldmath $ V$}_0^*, 
	\mbox{\scriptsize \boldmath $ V$}_0}
	+ a\bigl( \mbox{\boldmath $ \mu $}_\varepsilon (t),\mbox{\boldmath $ z$} \bigr) 
	=0 
	\quad \mbox{\it for all } \mbox{\boldmath $ z$} \in \mbox{\boldmath $ V$}_0,
	\label{weak1} 
	\\
	\bigl( 
	\mbox{\boldmath $ \mu $}_\varepsilon (t),\mbox{\boldmath $ z$}
	\bigr)_{\mbox{\scriptsize \boldmath $ H$}} 
	= \varepsilon a\bigl( \mbox{\boldmath $ v$}_\varepsilon (t),\mbox{\boldmath $ z$} \bigr) 
	+ \bigl( 
	 \mbox{\boldmath $ \xi $}_\varepsilon (t)
	+ \varepsilon \mbox{\boldmath $ \pi$} \bigl( 
	\mbox{\boldmath $ v$}_\varepsilon (t)+m_0\mbox{\boldmath $ 1$}
	\bigr)
	-\mbox{\boldmath $ f$}(t),\mbox{\boldmath $ z$}
	\bigr)_{\mbox{\scriptsize \boldmath $ H$}} 
	\nonumber 
	\\
	\quad \mbox{\it for all }\mbox{\boldmath $ z$} \in \mbox{\boldmath $ V$},
	\label{weak2}
\end{gather}
{\it for a.a.\ $t\in (0,T)$, with} 
\begin{equation*} 
	\xi _\varepsilon \in 
	\beta \bigl( v_\varepsilon +m_0 \bigr)
	\quad {\it a.e.~in~} Q,
	\quad 
	\xi _{\Gamma, \varepsilon} \in 
	\beta \bigl( v_{\Gamma, \varepsilon }+m_0 \bigr)
	\quad {\it a.e.~on~}\Sigma.
\end{equation*} 
{\it Moreover, 
$\mbox{\boldmath $ v$}_\varepsilon (0)=\mbox{\boldmath $ v$}_{0}$ 
in $\mbox{\boldmath $ H$}_0$ holds.} \\

We call the solution obtained by this proposition 
a kind of weak solution for the problem in \eqref{CHae1}--\eqref{CHae5}. 
This proposition is 
a direct consequence of \cite[Theorem~2.2]{CF15}. 
Indeed, 
we assumed (A1)--(A4), and then, 
from the construction \eqref{ax} of $\mbox{\boldmath $ f$}$, we see that 
$\mbox{\boldmath $ f$} \in L^2(0,T;\mbox{\boldmath $ V$}_0)$. 
Thus, all of the conditions needed 
to apply \cite[Theorem~2.2]{CF15} can be corrected.

To obtain uniform estimates with respect to $\varepsilon \in (0,1]$, 
we must consider the approximate problem for (P)$_\varepsilon $, 
which is the same strategy used to prove Proposition~2.1. 
Therefore, we give only a sketch of the proof here. 
For each $\lambda \in (0,1]$, consider the 
approximate problem (P)$_{\varepsilon,\lambda}$
\begin{align} 
	& \bigl \langle \mbox{\boldmath $ v$}'_{\varepsilon, \lambda } (t),\mbox{\boldmath $ z$} 
	\bigr \rangle _{\mbox{\scriptsize \boldmath $ V$}_0^*, 
	\mbox{\scriptsize \boldmath $ V$}_0}
	+ a\bigl( \mbox{\boldmath $ \mu $}_{\varepsilon,\lambda } (t),\mbox{\boldmath $ z$} \bigr) 
	=0 
	\quad \mbox{\rm for all } 
	\mbox{\boldmath $ z$} \in \mbox{\boldmath $ V$}_0,
	\nonumber 
	% \label{weakap1} 
	\\
	\bigl( 
	\mbox{\boldmath $ \mu $}_{\varepsilon,\lambda } (t),\mbox{\boldmath $ z$}
	\bigr)_{\mbox{\scriptsize \boldmath $ H$}} = {}
	& \lambda 
	\bigl(\mbox{\boldmath $ v$}_{\varepsilon, \lambda  }'(t),\mbox{\boldmath $ z$}
	\bigr)_{\mbox{\scriptsize \boldmath $ H$}}
	+ 
	\varepsilon a \bigl( 
	\mbox{\boldmath $ v$}_{\varepsilon,\lambda } (t),\mbox{\boldmath $ z$} \bigr) 
	\nonumber 
	\\
	& {} + \bigl( 
	\mbox{\boldmath $ \beta $}_\lambda  
	\bigl( \mbox{\boldmath $ v$}_{\varepsilon, \lambda } (t)+m_0\mbox{\boldmath $ 1$} \bigr) 
	+ \varepsilon \mbox{\boldmath $ \pi$} \bigl( 
	\mbox{\boldmath $ v$}_{\varepsilon,\lambda } (t)+m_0\mbox{\boldmath $ 1$}
	\bigr)
	-\mbox{\boldmath $ f$}(t),\mbox{\boldmath $ z$}
	\bigr)_{\mbox{\scriptsize \boldmath $ H$}} 
	\quad \mbox{\rm for all }\mbox{\boldmath $ z$} \in \mbox{\boldmath $ V$},
	\label{weakap2}
\end{align}
with $\mbox{\boldmath $ v$}_{\varepsilon,\lambda}(0) =\mbox{\boldmath $ v$}_0$ 
in $\mbox{\boldmath $ H$}_0$, 
where 
$\mbox{\boldmath $ \beta  $}_\lambda (\mbox{\boldmath $ z $})
:= (\beta _\lambda (z),\beta _\lambda (z_\Gamma ))$ 
for all $\mbox{\boldmath $ z $} \in \mbox{\boldmath $ H $}$ and 
$\beta _\lambda $ is the {Y}osida 
approximation 
(see, e.g., \cite{Bar10, Bre73, Ken07}). Namely, 
$\beta _\lambda  :\mathbb{R} \to \mathbb{R}$ 
along with the associated resolvent operator 
$J_\lambda :\mathbb{R} \to \mathbb{R}$ are defined 
by 
\begin{equation*}
	\beta _\lambda  (r)
	:= \frac{1}{\lambda  } \bigl( r-J_\lambda (r) \bigr)
	:=\frac{1}{\lambda  }\bigl( r-(I+ \lambda  \beta )^{-1} (r) \bigr)
	\quad {\rm for~all~} r \in \mathbb{R}. 
\end{equation*}
Moreover, the related {M}oreau--{Y}osida 
regularization $\widehat{\beta }_\lambda$
of $\widehat{\beta }:\mathbb{R} \to \mathbb{R}$ fulfills
\begin{equation}
	\widehat{\beta }_{\lambda }(r)
	:=\inf_{l \in \mathbb{R}}
	\left\{ \frac{1}{2 \lambda} |r-l|^2
	+\widehat{\beta }(l) \right\} 
	= 
	\frac{1}{2 \lambda } 
	\bigl| r-J_\lambda (r) \bigr|^2 
	+ \widehat{\beta }\bigl (J_\lambda (r) \bigr )
	= \int_{0}^{r} \beta _\lambda (l)dl
	\label{resol}
\end{equation}
for all $r \in \mathbb{R}$. 
We know the basic property
\begin{equation}
	\quad 0 \le \widehat{\beta }_\lambda (r) 
	\le \widehat{\beta }(r)
	\quad \mbox{for all } r \in \mathbb{R}.
	\label{prim}
\end{equation}

The problem (P)$_{\varepsilon ,\lambda }$ can be solved (see, e.g., \cite{CF14, CF15, CV90, KN96, Kub12}); 
more precisely, there exist
$\mbox{\boldmath $ v$}_{\varepsilon ,\lambda } 
\in H^1(0,T;\mbox{\boldmath $ H$}_0) \cap L^\infty (0,T;\mbox{\boldmath $ V$}_0)
\cap L^2(0,T;\mbox{\boldmath $ W$})$ and $\mbox{\boldmath $ \mu $}_{\varepsilon ,\lambda } \in 
L^2(0,T;\mbox{\boldmath $ V$})$ 
that satisfy
\begin{align} 
	\lambda \mbox{\boldmath $ v$}'_{\varepsilon ,\lambda }(t) 
	& + \mbox{\boldmath $ F$}^{-1}\mbox{\boldmath $ v$}'_{\varepsilon ,\lambda }(t)
	+ \varepsilon \partial \varphi \bigl( 
	\mbox{\boldmath $ v$}_{\varepsilon ,\lambda } (t) \bigr) 
	\nonumber \\
	& 
	=\mbox{\boldmath $ P$} \bigl( 
	-\mbox{\boldmath $ \beta $}_\lambda \bigl( 
	\mbox{\boldmath $ v$}_{\varepsilon ,\lambda }(t) + m_0 \mbox{\boldmath $ 1$} \bigr)
	- \varepsilon \mbox{\boldmath $ \pi $} 
	\bigl( \mbox{\boldmath $ v$}_{\varepsilon ,\lambda } (t) + m_0 \mbox{\boldmath $ 1$} \bigr) +\mbox{\boldmath $ f$}(t) \bigr)
	\quad {\rm in}~\mbox{\boldmath $ H$}_0
	\label{ee}
\end{align} 
for a.a.\ $t\in (0,T)$ 
with $\mbox{\boldmath $ v$}_{\varepsilon,\lambda } (0)=\mbox{\boldmath $ v$}_{0}$ in 
$\mbox{\boldmath $ H$}_0$, where 
$\varphi: \mbox{\boldmath $ H$}_0 \to [0,+\infty ]$ is a proper, 
lower semicontinuous, and convex functional 
\begin{equation}
	\varphi (\mbox{\boldmath $ z$}) 
	:= \begin{cases}
	\displaystyle 
	\frac{1}{2} \int_{\Omega }^{} | \nabla z |^2 dx 
	+\frac{1}{2} \int_{\Gamma }^{} 
	 |\nabla _{\Gamma } z_{\Gamma }  |^2 d\Gamma  
	\quad 
	\mbox{if } 
	\mbox{\boldmath $ z$} \in \mbox{\boldmath $ V$}_0, \vspace{2mm}\\
	+\infty \quad \mbox{otherwise}. 
	\end{cases} 
	\label{phi}
\end{equation}
Here, the subdifferential 
$\partial \varphi $ on $\mbox{\boldmath $ H$}_0$ is characterized by 
$\partial \varphi (\mbox{\boldmath $ z$})
=(-\Delta z,\partial_{\mbox{\scriptsize \boldmath $ \nu $}} z-\Delta _\Gamma z_\Gamma ) 
$ with 
$\mbox{\boldmath $ z$} \in 
D(\partial \varphi )=\mbox{\boldmath $ W$} \cap \mbox{\boldmath $ V$}_0$ 
(see, e.g., \cite[Lemma C]{CF15}). 
We also note that 
\begin{equation} 
	2 \varphi (\mbox{\boldmath $ z$})=  a(\mbox{\boldmath $ z$}, \mbox{\boldmath $ z$}  ) =
	\langle \mbox{\boldmath $ F$} \mbox{\boldmath $ z$}, \mbox{\boldmath $ z$} 
	\rangle _{\mbox{\scriptsize \boldmath $ V$}^*_0, 
	\mbox{\scriptsize \boldmath $ V$}_0}=|\mbox{\boldmath $ z$}|_{\mbox{\scriptsize \boldmath $ V$}_0}^2
	\quad \mbox{for all }\mbox{\boldmath $ z$} \in \mbox{\boldmath $ V$}_0.
	\label{star}
\end{equation}
Moreover, $\mbox{\boldmath $ P$}:\mbox{\boldmath $ H$} \to \mbox{\boldmath $ H$}_0$ is a projection defined by 
$\mbox{\boldmath $ P$}\mbox{\boldmath $ z$}:=\mbox{\boldmath $ z$}-m(\mbox{\boldmath $ z$}) \mbox{\boldmath $ 1$}$ for all 
$\mbox{\boldmath $ z$} \in \mbox{\boldmath $ H$}$, and it satisfies 
\begin{align*}
	(\mbox{\boldmath $ z$}_0,\mbox{\boldmath $ P$} \tilde{\mbox{\boldmath $ z$}})_{\mbox{\scriptsize \boldmath $ H$}_0}
	= (\mbox{\boldmath $ z$}_0,\tilde{\mbox{\boldmath $ z$}}
	)_{\mbox{\scriptsize \boldmath $ H$}}
	& \quad {\rm for~all~}\mbox{\boldmath $ z$}_0 \in \mbox{\boldmath $ H$}_0
	\ {\rm and}~\tilde{\mbox{\boldmath $ z$}} \in \mbox{\boldmath $ H$},\\
	|\mbox{\boldmath $ P$} \mbox{\boldmath $ z$}
	|_{\mbox{\scriptsize \boldmath $ V$}_0}
	= 
	|\mbox{\boldmath $ z$}|_{\mbox{\scriptsize \boldmath $ V$}_0}
	\le 
	|\mbox{\boldmath $ z$}|_{\mbox{\scriptsize \boldmath $ V$}}
	& \quad {\rm for~all~}\mbox{\boldmath $ z$} \in \mbox{\boldmath $ V$}.
\end{align*} 
The 
standard strategy obtains a priori estimates with respect 
to $\lambda \in (0,1]$ 
and considers the limiting procedure $\lambda \to 0$.

\paragraph{Remark 3.} 
Let $\mbox{\boldmath $ u$}_\varepsilon :=\mbox{\boldmath $ v$}_\varepsilon +m_0 \mbox{\boldmath $ 1$}$. 
Then, \cite[Remark 2]{CF15} means that \eqref{weak1} implies 
\begin{equation*} 
	\bigl \langle u'_\varepsilon (t),z 
	\bigr \rangle _{V^*,V}
	+ \langle u'_{\Gamma, \varepsilon} (t),z_\Gamma 
	\bigr \rangle _{V_\Gamma ^*,V_\Gamma }
	+ \int_{\Omega }^{} \nabla \mu _\varepsilon (t) \cdot \nabla z dx 
	+ \int_{\Gamma }^{} \nabla_\Gamma  \mu _{\Gamma, \varepsilon }(t) 
	\cdot \nabla _\Gamma z_\Gamma 
	d \Gamma 
	=0
%	\label{weakform1u}
\end{equation*}
for all $\mbox{\boldmath $ z$} \in \mbox{\boldmath $ V$}$. 
Moreover, thanks to the regularity 
$\mbox{\boldmath $ v$}_\varepsilon  \in L^2(0,T;\mbox{\boldmath $ W$})$, 
we see that \eqref{weak2} implies 
\begin{gather*} 
	\xi_\varepsilon  \in \beta (u_\varepsilon ), 
	\quad \mu_\varepsilon  = - \varepsilon \Delta u_\varepsilon  
	+ \xi_\varepsilon + \varepsilon \pi (u_\varepsilon )-f 
	\quad \mbox{a.e.\ in } Q,
	\\
	\xi_{\Gamma, \varepsilon}  \in \beta (u_{\Gamma, \varepsilon} ), 
	\quad 
	\mu _{\Gamma,\varepsilon  } = 
	\varepsilon \partial_{\mbox{\scriptsize \boldmath $ \nu $}} u_\varepsilon 
	- \varepsilon \Delta _\Gamma u_{\Gamma,\varepsilon } 
	+ \xi_{\Gamma,\varepsilon } 
	+ \varepsilon \pi _\Gamma (u_{\Gamma,\varepsilon}) 
	- f_\Gamma 
	\quad \mbox{a.e.\ on } \Sigma.
\end{gather*}

%%%%% Section 3. %%%%%
\section{Uniform estimates}
\setcounter{equation}{0}

In this section, we 
obtain the uniform estimates needed to prove 
the main theorem.

%%%%% Section 3.1. %%%%%
\subsection{Uniform estimates for approximate solutions of (P)$_{\varepsilon,\lambda}$}

To prove the main theorem, we will use the uniform estimates that are independent 
of $\varepsilon $ for the solutions of (P)$_\varepsilon $. 
For this, we must start from 
the approximate problem (P)$_{\varepsilon,\lambda}$. 
In this subsection, we obtain uniform estimates 
for the approximate solutions of (P)$_{\varepsilon,\lambda}$.
We recall the change of variable $\mbox{\boldmath $ u$}_{\varepsilon ,\lambda }:=
\mbox{\boldmath $ v$}_{\varepsilon ,\lambda }+m_0 \mbox{\boldmath $ 1$}$.

\paragraph{Lemma 3.1.}
{\it There exist positive constants $M_1$ and $M_2$, 
independent of $\varepsilon \in (0,1/4]$ and $\lambda \in (0,1]$, such that}
\begin{gather*}
	\lambda 
	\bigl| \mbox{\boldmath $ v$}_{\varepsilon,\lambda}(t)
	\bigr|_{\mbox{\scriptsize \boldmath $ H$}_0}^2 
	+ 
	\bigl| \mbox{\boldmath $ v$}_{\varepsilon,\lambda }(t)
	\bigr|_{\mbox{\scriptsize \boldmath $ V$}_0^*}^2
	\le M_1, \\
	\frac{\varepsilon}{2} \int_{0}^{t} 
	\bigl |\mbox{\boldmath $ v$}_{\varepsilon,\lambda  }(s)
	\bigr |_{\mbox{\scriptsize \boldmath $ V$}_0}^2 ds 
	+ 2 \int_{0}^{t} \bigl| \widehat{\beta}_\lambda  \bigl(u_{\varepsilon,\lambda }(s) \bigr) 
	\bigr|_{L^1(\Omega )} ds 
	+ 2 \int_{0}^{t} \bigl| 
	\widehat{\beta} _\lambda \bigl( u_{\Gamma, \varepsilon, \lambda}(s) \bigr) 
	\bigr|_{L^1(\Gamma)} ds
	\le M_2
\end{gather*}
{\it for all $t\in [0,T]$.}

\paragraph{Proof.} We test \eqref{ee} at time $s \in (0,T)$ by 
$\mbox{\boldmath $ v$}_{\varepsilon,\lambda} (s) \in \mbox{\boldmath $ V$}_0$ 
which is considered in problem (P)$_{\varepsilon ,\lambda }$.  
Then, using \eqref{Vstar}, we have 
\begin{align}
	& \lambda 
	\bigl( 
	\mbox{\boldmath $ v$}_{\varepsilon ,\lambda }' (s),
	\mbox{\boldmath $ v$}_{\varepsilon ,\lambda } (s) 
	\bigr)_{\mbox{\scriptsize \boldmath $ H$}_0}  
	+  \bigl( 
	\mbox{\boldmath $ v$}_{\varepsilon ,\lambda }' (s),
	\mbox{\boldmath $ v$}_{\varepsilon ,\lambda } (s)
	\bigr )_{\mbox{\scriptsize \boldmath $ V$}_0^*}
	+ \varepsilon \bigl( \partial \varphi \bigl ( \mbox{\boldmath $ v$}_{\varepsilon,\lambda } (s) \bigr ), 
	\mbox{\boldmath $ v$}_{\varepsilon ,\lambda }(s) \bigr)_{\mbox{\scriptsize \boldmath $ H$}_0} 
	\nonumber \\
	& \quad {} + \bigl( \mbox{\boldmath $ P$} \mbox{\boldmath $ \beta $}_\lambda
	\bigl ( \mbox{\boldmath $ v$}_{\varepsilon,\lambda } (s)+m_0\mbox{\boldmath $ 1$} \bigr ), 
	\mbox{\boldmath $ v$}_{\varepsilon,\lambda } (s) 
	\bigr)_{\mbox{\scriptsize \boldmath $ H$}_0} 
	= 
	\bigl(\mbox{\boldmath $ f$}(s)-\varepsilon \mbox{\boldmath $ P$}
	\mbox{\boldmath $ \pi $}
	\bigl( \mbox{\boldmath $ v$}_{\varepsilon ,\lambda } (s) \bigr), 
	\mbox{\boldmath $ v$}_{\varepsilon,\lambda } (s) 
	\bigr)_{\mbox{\scriptsize \boldmath $ H$}_0}
	\label{1st} 
\end{align} 
for a.a.\ $s \in (0,T)$. 
Now, from the definition of the subdifferential, we have 
\begin{align}
	& 
	\bigl( \mbox{\boldmath $ P$} \mbox{\boldmath $ \beta $}_\lambda  
	\bigl ( \mbox{\boldmath $ v$}_{\varepsilon,\lambda } (s)+m_0\mbox{\boldmath $ 1$} \bigr), 
	\mbox{\boldmath $ v$}_{\varepsilon,\lambda }(s) 
	\bigr)_{\mbox{\scriptsize \boldmath $ H$}_0} 
	\nonumber \\
	& =  \int_{\Omega }^{} \beta_\lambda  
	\bigl( u_{\varepsilon,\lambda}(s) \bigr) 
	\bigl( u_{\varepsilon,\lambda } (s) - m_0 \bigr) dx 
	+ \int_{\Gamma}^{} \beta _\lambda 
	\bigl( u_{\Gamma, \varepsilon, \lambda } (s) \bigr) 
	\bigl( u_{\Gamma, \varepsilon, \lambda } (s) - m_0 \bigr) d\Gamma 
	\nonumber \\
	& \ge 
	\int_{\Omega }^{} \widehat{\beta} _\lambda 
	\bigl ( u_{\varepsilon,\lambda}(s) \bigr ) dx 
	 - 
	 \widehat{\beta }_\lambda 
	(m_0) |\Omega |
	+ \int_{\Gamma}^{} \widehat{\beta}_\lambda  
	\bigl ( u_{\Gamma, \varepsilon, \lambda } (s) \bigr ) d\Gamma 
	-  \widehat{\beta }_\lambda 
	(m_0) |\Gamma |
	\label{lem1}
\end{align}
for a.a.\ $s \in (0,T)$.
Next, recalling the fundamental property of the chain rule with \eqref{pi} 
and using the Young inequality, we have 
\begin{align}
	& - \varepsilon \bigl \langle  
	\mbox{\boldmath $ v$}_{\varepsilon,\lambda }(s), 
	 \mbox{\boldmath $ P$} \mbox{\boldmath $ \pi $} 
	\bigl ( 
	\mbox{\boldmath $ v$}_{\varepsilon,\lambda } (s)
	+m_0\mbox{\boldmath $ 1$} 
	\bigr) 
	\bigr \rangle _{\mbox{\scriptsize \boldmath $ V$}_0^*,\mbox{\scriptsize \boldmath $ V$}_0}
	\nonumber \\
	& \le 
	\frac{1}{4} 
	\bigl| 
	\mbox{\boldmath $ v$}_{\varepsilon,\lambda }(s)
	\bigr|_{\mbox{\scriptsize \boldmath $ V$}_0^*}^2 
	+ \varepsilon ^2 
	 \int_{\Omega }^{} \bigl| \nabla \pi 
	\bigl( v_{\varepsilon,\lambda}(s) +m_0 \bigr) \bigr|^2 dx 
	+ \varepsilon ^2 \int_{\Gamma}^{} \bigl| \nabla _\Gamma \pi 
	\bigl( v_{\Gamma, \varepsilon, \lambda } (s)+m_0 \bigr) \bigr|^2 d\Gamma 
	\nonumber \\
	& \le 
	\frac{1}{4} 
	\bigl| 
	\mbox{\boldmath $ v$}_{\varepsilon,\lambda }(s)
	\bigr|_{\mbox{\scriptsize \boldmath $ V$}_0^*}^2 
	+ \varepsilon ^2 \bigl| \mbox{\boldmath $ v$}_{\varepsilon ,\lambda }(s) \bigr|_{\mbox{\scriptsize \boldmath $ V$}_0}^2,  
	\label{lem2}
\end{align}
and 
\begin{equation} 
	\bigl \langle  
	\mbox{\boldmath $ v$}_{\varepsilon,\lambda }(s), 
	  \mbox{\boldmath $ f$}(s)
	\bigr \rangle _{\mbox{\scriptsize \boldmath $ V$}_0^*,\mbox{\scriptsize \boldmath $ V$}_0}
	\le 
	\frac{1}{4} 
	\bigl| 
	\mbox{\boldmath $ v$}_{\varepsilon,\lambda }(s)
	\bigr|_{\mbox{\scriptsize \boldmath $ V$}_0^*}^2 
	+ \bigl| \mbox{\boldmath $ f$}(s) \bigr|_{\mbox{\scriptsize \boldmath $ V$}_0}^2  
	\label{lem3}
\end{equation}
for a.a.\ $s \in (0,T)$.
Thus, correcting \eqref{1st}--\eqref{lem3},
recalling  
the definition of the subdifferential and 
using \eqref{prim}, \eqref{star}, and 
$\varphi(\mbox{\boldmath $ 0$})=0$, 
we see that 
\begin{align}
	& \lambda 
	\frac{d}{ds}
	\bigl | 
	\mbox{\boldmath $ v$}_{\varepsilon, \lambda } (s )
	\bigr|_{\mbox{\scriptsize \boldmath $ H$}_0}^2 
	+ 
	\frac{d}{ds} 
	\bigl | 
	\mbox{\boldmath $ v$}_{\varepsilon, \lambda }(s )
	\bigr |_{\mbox{\scriptsize \boldmath $ V$}_0^*}^2
	+ (\varepsilon -2\varepsilon ^2) 
	\bigl | \mbox{\boldmath $ v$}_{\varepsilon,\lambda }(s )
	\bigr |_{\mbox{\scriptsize \boldmath $ V$}_0}^2
	\nonumber \\
	&	
	\quad {} + 2
	\bigl| \widehat{\beta}_\lambda 
	\bigl ( u_{\varepsilon,\lambda } (s) 
	\bigr) \bigr |_{L^1(\Omega )}
	+ 2 \bigl| \widehat{\beta}_\lambda 
	\bigl ( u_{\Gamma, \varepsilon,\lambda } (s) 
	\bigr) \bigr |_{L^1(\Gamma )}
	\nonumber 
	\\
	& \le 
	2\bigl( |\Omega |+|\Gamma | \bigr) \widehat{\beta }(m_0)
	+ \bigl| 
	\mbox{\boldmath $ v$}_{\varepsilon,\lambda }(s)
	\bigr|_{\mbox{\scriptsize \boldmath $ V$}_0^*}^2 
	+ 2
	\bigl| \mbox{\boldmath $ f$}(s) \bigr|_{\mbox{\scriptsize \boldmath $ V$}_0}^2
	\label{lem4}
\end{align}
for a.a.\ $s \in (0,T)$. 
If we take $\varepsilon \in (0,1/4]$, then 
$\varepsilon -2\varepsilon ^2 \ge \varepsilon /2>0$. 
By virtue of the Gronwall inequality, we obtain
\begin{align*}
	& \lambda 
	\bigl | 
	\mbox{\boldmath $ v$}_{\varepsilon, \lambda }(t)
	\bigr|_{\mbox{\scriptsize \boldmath $ H$}_0}^2 
	+  
	\bigl | 
	\mbox{\boldmath $ v$}_{\varepsilon, \lambda }(t)
	\bigr |_{\mbox{\scriptsize \boldmath $ V$}_0^*}^2
	\nonumber \\
	& \le 
	\left( 
	| 
	\mbox{\boldmath $ v$}_0
	|_{\mbox{\scriptsize \boldmath $ H$}_0}^2 
	+  
	| 
	\mbox{\boldmath $ v$}_0
	|_{\mbox{\scriptsize \boldmath $ V$}_0^*}^2
	+
	2 T \bigl( |\Omega |+|\Gamma | \bigr) \widehat{\beta }(m_0)
	+ 2 
	| \mbox{\boldmath $ f$} 
	|_{L^2(0,T;\mbox{\scriptsize \boldmath $ V$}_0)}^2 
	\right) e^T =:M_1
\end{align*} 
for all $t \in [0,T]$. 
Next, integrating \eqref{lem4} over $(0,t)$ with respect 
to $s$, we obtain
\begin{align*}
	& 
	\frac{\varepsilon }{2}
	\int_{0}^{t} 
	\bigl | \mbox{\boldmath $ v$}_{\varepsilon,\lambda }(s )
	\bigr |_{\mbox{\scriptsize \boldmath $ V$}_0}^2ds 
	+ 2
	\int_{0}^{t} \bigl| \widehat{\beta} _\lambda 
	\bigl ( u_{\varepsilon,\lambda } (s) 
	\bigr) \bigr |_{L^1(\Omega )}ds 
	+ 2 \int_{0}^{t} 
	\bigl| \widehat{\beta}_\lambda 
	\bigl ( u_{\Gamma, \varepsilon,\lambda } (s) 
	\bigr) \bigr |_{L^1(\Gamma )}ds
	\nonumber 
	\\
	& \le  
	|\mbox{\boldmath $ v $}_0|_{\mbox{\scriptsize \boldmath $ H $}_0}^2
	+ 
	|\mbox{\boldmath $ v $}_0|_{\mbox{\scriptsize \boldmath $ V $}_0^*}^2
	+
	2T\bigl( |\Omega |+|\Gamma | \bigr) \widehat{\beta }(m_0)
	+ 2
	| \mbox{\boldmath $ f$}
	|_{L^2(0,T;\mbox{\scriptsize \boldmath $ V$}_0)}^2  
	+ 
	\int_{0}^{t} 
	\bigl| 
	\mbox{\boldmath $ v$}_{\varepsilon,\lambda }(s)
	\bigr|_{\mbox{\scriptsize \boldmath $ V$}_0^*}^2 
	ds
	\nonumber \\
	& \le 
	M_1 (1+T):=M_2
\end{align*} 
for all $t \in [0,T]$. Thus, we have the conclusion. \hfill $\Box$

\paragraph{Lemma 3.2.}
{\it There exists a positive constant $M_3$, 
independent of $\varepsilon \in (0,1/4]$ and $\lambda \in (0,1]$, such that}
\begin{align}
	& 
	2 \lambda 
	\int_{0}^{t} 
	\bigl| 
	\mbox{\boldmath $ v$}_{\varepsilon,\lambda } '(s) 
	\bigr|_{\mbox{\scriptsize \boldmath $ H$}_0}^2 
	ds
	+ 
	\int_{0}^{t} 
	\bigl| 
	\mbox{\boldmath $ v$}_{\varepsilon,\lambda } '(s)
	| _{\mbox{\scriptsize \boldmath $ V$}_0^*}^2 ds 
	+ 
	\varepsilon 
	\bigl|\mbox{\boldmath $ v$}_{\varepsilon,\lambda }(t) 
	\bigr|_{\mbox{\scriptsize \boldmath $ V$}_0}^2  
	\nonumber
	\\ 
	&{}+ 2 \bigl| \widehat{\beta}_\lambda  \bigl (u_{\varepsilon,\lambda}(t) \bigr ) 
	\bigr|_{L^1(\Omega )} 
	+ 2 \bigl| \widehat{\beta}_\lambda  \bigl ( 
	u_{\Gamma, \varepsilon,\lambda }(t) \bigr) 
	\bigr|_{L^1(\Gamma )} \le M_3, 
	\label{m3a}
\end{align}
\begin{equation} 
	 \int_{0}^{t} \bigl| \mbox{\boldmath $ P$} 
	 \mbox{\boldmath $ \mu $}_{\varepsilon, \lambda }(s) 
	 \bigr|_{\mbox{\scriptsize \boldmath $ V$}_0}^2 ds 
	 \le M_3
	 \label{m3b}
\end{equation} 
{\it for all $t \in [0,T]$.}
	 
\paragraph{Proof.} 
We test \eqref{ee} at time $s \in (0,T)$ by 
$\mbox{\boldmath $ v$}_{\varepsilon,\lambda }' (s) 
\in \mbox{\boldmath $ H$}_0$. 
Then, using the same method as for \eqref{lem2}, we have 
\begin{equation*}
	- \varepsilon \bigl \langle  
	\mbox{\boldmath $ v$}_{\varepsilon,\lambda }'(s), 
	 \mbox{\boldmath $ P$} \mbox{\boldmath $ \pi $} 
	\bigl ( 
	\mbox{\boldmath $ v$}_{\varepsilon,\lambda } (s)
	+m_0\mbox{\boldmath $ 1$} 
	\bigr) 
	\bigr \rangle _{\mbox{\scriptsize \boldmath $ V$}_0^*,\mbox{\scriptsize \boldmath $ V$}_0}
	\le 
	\frac{1}{4} 
	\bigl| 
	\mbox{\boldmath $ v$}_{\varepsilon,\lambda }'(s)
	\bigr|_{\mbox{\scriptsize \boldmath $ V$}_0^*}^2 
	+ \varepsilon ^2 \bigl| \mbox{\boldmath $ v$}_{\varepsilon ,\lambda }(s) \bigr|_{\mbox{\scriptsize \boldmath $ V$}_0}^2,  
\end{equation*}
and 
\begin{equation*} 
	\bigl \langle  
	\mbox{\boldmath $ v$}_{\varepsilon,\lambda }'(s), 
	\mbox{\boldmath $ f$}
	(s)
	\bigr \rangle _{\mbox{\scriptsize \boldmath $ V$}_0^*,\mbox{\scriptsize \boldmath $ V$}_0}
	\le 
	\frac{1}{4} 
	\bigl| 
	\mbox{\boldmath $ v$}_{\varepsilon,\lambda }'(s)
	\bigr|_{\mbox{\scriptsize \boldmath $ V$}_0^*}^2 
	+ \bigl| \mbox{\boldmath $ f$}(s) \bigr|_{\mbox{\scriptsize \boldmath $ V$}_0}^2  
\end{equation*}
for a.a.\ $s \in (0,T)$. 
Then, using \eqref{Vstar} and the chain rule, we deduce
\begin{align}
	& \lambda \bigl 
	|\mbox{\boldmath $ v$}_{\varepsilon,\lambda } '(s ) 
	\bigr |_{\mbox{\scriptsize \boldmath $ H$}_0}^2  
	+  
	\frac{1}{2} \bigl |\mbox{\boldmath $ v$}_{\varepsilon,\lambda } '(s) 
	\bigr |_{\mbox{\scriptsize \boldmath $ V$}_0^*}^2 
	+ \varepsilon \frac{d}{ds} 
	\varphi \bigl ( \mbox{\boldmath $ v$}_{\varepsilon,\lambda } (s) \bigr )
	\nonumber \\
	&\quad {}+\frac{d}{ds} \int_{\Omega }^{} 
	\widehat{\beta }_\lambda  \bigl (
	v_{\varepsilon,\lambda } (s)+m_0)dx
	+
	\frac{d}{ds} \int_{\Gamma }^{} 
	\widehat{\beta }_\lambda  
	\bigl (v_{\Gamma, \varepsilon,\lambda } (s)+m_0\bigr )d\Gamma 
	\nonumber \\
	& \le 
	\varepsilon ^2 \bigl| \mbox{\boldmath $ v$}_{\varepsilon ,\lambda }(s) \bigr|_{\mbox{\scriptsize \boldmath $ V$}_0}^2
	+ \bigl| \mbox{\boldmath $ f$}(s) \bigr|_{\mbox{\scriptsize \boldmath $ V$}_0}^2  
	\label{lem5}
\end{align}
for a.a.\ $s \in (0,T)$. 
Integrating \eqref{lem5} over $(0,t)$ with respect to $s$, 
we can find a positive constant $M_3$, depending only on 
$|\mbox{\boldmath $ v$}_0|_{\mbox{\scriptsize \boldmath $ V$}_0}$,  
$|\widehat{\beta }(u_0)|_{L^1(\Omega )}$, 
$|\widehat{\beta }(u_{0\Gamma })|_{L^1(\Gamma )}$, 
$M_2$, and 
$|\mbox{\boldmath $ f$}|_{L^2(0,T;\mbox{\scriptsize \boldmath $ V$}_0)}$,
such that the aforementioned estimate \eqref{m3a} holds. 
Next, to obtain \eqref{m3b}, 
we recall \eqref{weakap2}, (see also Remark~3). 
We have 
\begin{equation}
	\mbox{\boldmath $ \mu $}_{\varepsilon,\lambda } (s)= 
	 \lambda 
	\mbox{\boldmath $ v$}_{\varepsilon, \lambda  }'(s)
	+ \varepsilon 
	\partial \varphi 
	\bigl( 
	\mbox{\boldmath $ v$}_{\varepsilon,\lambda } (s) \bigr) 
	+ 
	\mbox{\boldmath $ \beta $}_\lambda  
	\bigl( \mbox{\boldmath $ u$}_{\varepsilon, \lambda } (s) \bigr) 
	+ \varepsilon \mbox{\boldmath $ \pi$} \bigl( 
	\mbox{\boldmath $ u$}_{\varepsilon,\lambda } (s)
	\bigr)
	-\mbox{\boldmath $ f$}(s) 
	\quad {\rm in~} \mbox{\boldmath $ V$},
	\label{mu}
\end{equation}
for a.a.\ $s \in (0,T)$.  
Therefore, the evolution equation \eqref{ee} is equivalent to 
\begin{equation} 
	\mbox{\boldmath $ v$}_{\varepsilon ,\lambda }'(s)
	+ \mbox{\boldmath $ F$} \bigl( \mbox{\boldmath $ P$} \mbox{\boldmath $ \mu $}_{\varepsilon, \lambda }(s) 
	\bigr) = \mbox{\boldmath $ 0$}
	\quad {\rm in}~\mbox{\boldmath $ V$}_0^*, \quad {\rm for~a.a.~} s \in (0,T),
	\label{ee2}
\end{equation} 
with \eqref{mu}.
We test \eqref{ee2} by 
$\mbox{\boldmath $ P$} \mbox{\boldmath $ \mu $}_{\varepsilon,\lambda }(s) 
\in \mbox{\boldmath $ V$}_0$, and integrate the 
resultant over $(0,t)$ with respect to $s$. 
Then, using \eqref{dual} and the Young inequality, we have 
\begin{equation*} 
	\int_{0}^{t} \bigl| \mbox{\boldmath $ P$} \mbox{\boldmath $ \mu $}_{\varepsilon ,\lambda } (s) 
	\bigr|_{\mbox{\scriptsize \boldmath $ V$}_0}^2 ds 
	\le \frac{1}{2}
	\int_{0}^{t} 
	\bigl|
	\mbox{\boldmath $ v$}_{\varepsilon ,\lambda }'(s) 
	\bigr|_{\mbox{\scriptsize \boldmath $ V$}_0^*}^2 ds
	+ 
	\frac{1}{2}
	\int_{0}^{t} 
	\bigl| 
	\mbox{\boldmath $ P$} 
	\mbox{\boldmath $ \mu $}_{\varepsilon ,\lambda }(s) 
	\bigr|_{\mbox{\scriptsize \boldmath $ V$}_0}^2 ds
\end{equation*} 
for all $t \in [0,T]$. Thus, using \eqref{m3a}, we obtain \eqref{m3b}. 
\hfill $\Box$

\paragraph{Lemma 3.3.}
{\it There exist two positive constants $M_4$ and $M_5$,
independent of $\varepsilon \in (0,1/4]$ and $\lambda \in (0,1]$, such that}
\begin{gather}
	 \bigl| 
	 \mbox{\boldmath $ u $}_{\varepsilon, \lambda }(t) 
	 \bigr|_{\mbox{\scriptsize \boldmath $ H$}}^2
	 \le M_4,
	 \quad \bigl| 
	 \mbox{\boldmath $ v $}_{\varepsilon, \lambda }(t) 
	 \bigr|_{\mbox{\scriptsize \boldmath $ H$}_0}^2
	 \le M_4,
	 \label{m4a}
	\\
	\int_{0}^{t} 
	\bigl| \beta _\lambda \bigl (u_{\varepsilon,\lambda }(s) \bigr) 
	\bigr|_{L^1(\Omega )}^2 ds 
	+ 
	\int_{0}^{t} 
	\bigl| \beta _\lambda \bigl(u_{\Gamma, \varepsilon,\lambda} (s) \bigr) 
	\bigr|_{L^1(\Gamma)}^2 ds
	 \le M_5 \label{m4b}
\end{gather}
{\it for all $t\in [0,T]$.}

\paragraph{Proof.}
From \eqref{prim0} in (A1) and \eqref{resol}, 
\begin{align*}
	\widehat{\beta }_\lambda (r) 
	& = \frac{1}{2\lambda } 
	\bigl| r- J_\lambda (r) \bigr|^2 + \widehat{\beta } 
	\bigl( J_\lambda (r) \bigr) \\
	& \ge \frac{1}{2\lambda } 
	\bigl| r- J_\lambda (r) \bigr|^2 + 
	c_1 \bigl| J_\lambda (r) \bigr|^2 -c_2
	\quad {\rm for~all}~ r\in \mathbb{R}.
\end{align*}
Therefore, from \eqref{m3a} in Lemma~3.2, we obtain
\begin{gather*}
	\frac{1}{2\lambda } 
	\bigl| u_{\varepsilon ,\lambda }(t) - J_\lambda 
	\bigl( u_{\varepsilon ,\lambda }(t) 
	\bigr)  \bigr|_H^2 
	+c_1 \bigl| J_\lambda 
	\bigl( u_{\varepsilon ,\lambda }(t) 
	\bigr)  \bigr|_H^2 
	\le \frac{M_3}{2} +c_2 |\Omega |, \\
	\frac{1}{2\lambda } 
	\bigl| u_{\Gamma, \varepsilon ,\lambda }(t) -  J_\lambda 
	\bigl( u_{\Gamma, \varepsilon ,\lambda }(t) 
	\bigr) \bigr|_{H_\Gamma }^2 
	+c_1 \bigl| J_\lambda 
	\bigl( u_{\Gamma, \varepsilon ,\lambda }(t) 
	\bigr)  \bigr|_{H_\Gamma}^2 
	\le \frac{M_3}{2} +c_2 |\Gamma|
\end{gather*}
for all $t \in [0,T]$. Thanks to this $\lambda \in (0,1]$ implies that there 
exists $\tilde{M}_4>0$, depending only on 
$c_1$, $c_2$, $M_3$, $|\Omega |$, and $|\Gamma|$,
such that 
\begin{align*}
\bigl| \mbox{\boldmath $ u$}_{\varepsilon ,\lambda }(t) \bigr|_{\mbox{\scriptsize \boldmath $ H$}} ^2 
	& = 
	\bigl| u_{\varepsilon,\lambda}(t) \bigr|^2 _H
	+ 
	\bigl| u_{\Gamma, \varepsilon, \lambda } (t) \bigr|^2 _{H_\Gamma } \\
	& \le 
	2 \bigl| u_{\varepsilon ,\lambda }(t) - J_\lambda 
	\bigl( u_{\varepsilon ,\lambda }(t) 
	\bigr)  \bigr|_H^2 
	+2 \bigl| J_\lambda 
	\bigl( u_{\varepsilon ,\lambda }(t) 
	\bigr)  \bigr|_H^2 \nonumber \\
	& \quad {}
	+ 2 \bigl| u_{\Gamma, \varepsilon ,\lambda }(t) -  J_\lambda 
	\bigl( u_{\Gamma, \varepsilon ,\lambda }(t) 
	\bigr) \bigr|_{H_\Gamma }^2 
	+2 \bigl| J_\lambda 
	\bigl( u_{\Gamma, \varepsilon ,\lambda }(t) 
	\bigr)  \bigr|_{H_\Gamma}^2 
	\nonumber \\
	& \le \tilde{M}_4
\end{align*}
for all $t \in [0,T]$. Then, there exists $M_4>0$, depending only on 
$\tilde{M}_4$, $|m_0|$, $|\Omega |$, and $|\Gamma|$,
such that 
\begin{align*}
\bigl| \mbox{\boldmath $ v$}_{\varepsilon ,\lambda }(t) \bigr|_{\mbox{\scriptsize \boldmath $ H$}} ^2 
	& \le 2 \bigl| \mbox{\boldmath $ u$}_{\varepsilon ,\lambda }(t) \bigr|_{\mbox{\scriptsize \boldmath $ H$}} ^2 
	+ 2 |m_0|^2 \bigl( |\Omega |+|\Gamma |\bigr) \\
	& \le M_4
\end{align*}
for all $t \in [0,T]$. Thus, \eqref{m4a} has been obtained. 
Next, recall the useful inequality \cite[Section~5]{GMS09}. Indeed, from assumption 
(A4),
we have $m_0\in {\rm int}D(\beta )$ (which is the criterion needed to 
apply the inequality), meaning there exist two constants 
$c_3$, $c_4>0$ such that
\begin{equation} 
	\beta _\lambda (r)(r-m_0) \ge c_3 \beta _\lambda (r) -c_4
\end{equation} 
for all $r \in \mathbb{R}$ and $\lambda \in (0,1]$. 
Then, \eqref{lem1} can be improved by  
\begin{align*}
	& 
	\bigl( \mbox{\boldmath $ P$} \mbox{\boldmath $ \beta $}_\lambda  
	\bigl ( \mbox{\boldmath $ v$}_{\varepsilon,\lambda } (s)+m_0\mbox{\boldmath $ 1$} \bigr), 
	\mbox{\boldmath $ v$}_{\varepsilon,\lambda }(s) 
	\bigr)_{\mbox{\scriptsize \boldmath $ H$}_0} 
	\nonumber \\
	& =  \int_{\Omega }^{} \beta_\lambda  
	\bigl( u_{\varepsilon,\lambda}(s) \bigr) 
	\bigl( u_{\varepsilon,\lambda } (s) - m_0 \bigr) dx 
	+ \int_{\Gamma}^{} \beta _\lambda 
	\bigl( u_{\Gamma, \varepsilon, \lambda } (s) \bigr) 
	\bigl( u_{\Gamma, \varepsilon, \lambda } (s) - m_0 \bigr) d\Gamma 
	\nonumber \\
	& \ge c_3 
	\int_{\Omega }^{}  {\beta} _\lambda 
	\bigl ( u_{\varepsilon,\lambda}(s) \bigr)  dx 
	 - c_4 |\Omega |
	+ c_3 \int_{\Gamma}^{}  \beta_\lambda  
	\bigl ( u_{\Gamma, \varepsilon, \lambda } (s) \bigr)  d\Gamma 
	- c_4 |\Gamma |
\end{align*}
for a.a.\ $s \in (0,T)$.
Therefore, by using \eqref{Vstar} and the monotonicity of $\partial \varphi$,
\begin{align*}
	& c_3
	\int_{\Omega }^{} {\beta} _\lambda 
	\bigl ( u_{\varepsilon,\lambda}(s) \bigr)  dx 
	+ c_3 \int_{\Gamma}^{} \beta_\lambda  
	\bigl ( u_{\Gamma, \varepsilon, \lambda } (s) \bigr) d\Gamma 
	\nonumber \\
	& 
	\le c_4 \bigl(|\Omega | +  |\Gamma |\bigr)
	+ 
	\lambda 
	\bigl|
	\mbox{\boldmath $ v$}_{\varepsilon ,\lambda }' (s)
	\bigr|_{\mbox{\scriptsize \boldmath $ H$}_0} 
	\bigl| 
	\mbox{\boldmath $ v$}_{\varepsilon ,\lambda } (s) 
	\bigr|_{\mbox{\scriptsize \boldmath $ H$}_0}  
	+ 
	\bigl| 
	\mbox{\boldmath $ v$}_{\varepsilon ,\lambda }' (s)
	\bigr|_{\mbox{\scriptsize \boldmath $ V$}_0^*}
	\bigl| 
	\mbox{\boldmath $ v$}_{\varepsilon ,\lambda } (s)
	\bigr|_{\mbox{\scriptsize \boldmath $ V$}_0^*}
	\nonumber \\
	& \quad {}  
	+ \varepsilon 
	\bigl| 
	\mbox{\boldmath $ \pi $}
	\bigl( \mbox{\boldmath $ u$}_{\varepsilon ,\lambda } (s) 
	\bigr) 
	\bigr|_{\mbox{\scriptsize \boldmath $ H$}} 
	\bigl| 
	\mbox{\boldmath $ v$}_{\varepsilon,\lambda } (s) 
	\bigr|_{\mbox{\scriptsize \boldmath $ H$}_0}
	+ 
	\bigl| 
	\mbox{\boldmath $ f$}(s)
	\bigr| _{\mbox{\scriptsize \boldmath $ H$}_0}
	\bigl| 
	\mbox{\boldmath $ v$}_{\varepsilon,\lambda } (s) 
	\bigr|_{\mbox{\scriptsize \boldmath $ H$}_0}
\end{align*} 
for a.a.\ $s \in (0,T)$. 
Now, squaring the above and using Lemma~3.1 and \eqref{m4a}, we obtain
\begin{align}
	& \left( c_3 \int_{\Omega }^{}
	 \beta _\lambda \bigl( u_{\varepsilon,\lambda }(s) \bigr) 
	 dx
	+ 
	c_3 \int_{\Gamma }^{} 
	 \beta _\lambda  \bigl( u_{\Gamma, \varepsilon,\lambda }(s) \bigr) 
	 d\Gamma \right)^2
	 \nonumber \\
	& \le 5 c_4^{2} \bigl( |\Omega |+|\Gamma |\bigr)^2
	+ 5\lambda^2 M_4
	\bigl| 
	\mbox{\boldmath $ v$}_{\varepsilon ,\lambda }' (s)
	\bigr|_{\mbox{\scriptsize \boldmath $ H$}_0}^2
	+5 M_1 \bigl| 
	\mbox{\boldmath $ v $}_{\varepsilon,\lambda }' (s) 
	\bigr|_{\mbox{\scriptsize \boldmath $ V$}_0^*}^2
	\nonumber 
	\\
	& \quad {}
	+ 
	5 \varepsilon ^2 M_4
	\bigl| \mbox{\boldmath $ \pi $}
	\bigl (\mbox{\boldmath $ u$}_{\varepsilon, \lambda} (s) \bigr )
	\bigr|_{\mbox{\scriptsize \boldmath $ H$}}^2
	+
	5 M_4
	\bigl| \mbox{\boldmath $ f$}(s)
	\bigr| _{\mbox{\scriptsize \boldmath $ H$}_0}^2
	\label{sq}
\end{align}
for a.a.\ $s \in (0,T)$. 
Thus, from \eqref{pi} and \eqref{m4a}, we have 
\begin{align}
	\bigl| \mbox{\boldmath $ \pi $}
	\bigl (\mbox{\boldmath $ u$}_{\varepsilon, \lambda} (s) \bigr )
	\bigr|_{\mbox{\scriptsize \boldmath $ H$}}^2
	& = 
	\bigl| \mbox{\boldmath $ \pi $}
	\bigl (\mbox{\boldmath $ u$}_{\varepsilon, \lambda} (s) \bigr )
	-\mbox{\boldmath $ \pi $}
	(m_0\mbox{\boldmath $ 1$}) + \mbox{\boldmath $ \pi $}
	(m_0\mbox{\boldmath $ 1$})
	\bigr|_{\mbox{\scriptsize \boldmath $ H$}}^2
	\nonumber 
	\\
	& \le 2 \left\{  \int_{\Omega }^{}
	\Bigl( 
	\bigl| 
	v_{\varepsilon, \lambda} (s)
	\bigr|^2 
	+ 
	\bigl| 
	\pi (m_0)
	\bigr|^2
	\Bigr) 
	dx 
	+ \int_{\Gamma }^{}
	\Bigl( 
	\bigl| 
	v_{\Gamma, \varepsilon, \lambda} (s)
	\bigr|^2 
	+ 
	\bigl| 
	\pi (m_0)
	\bigr|^2
	\Bigr)  d\Gamma  
	\right\} 
	\nonumber \\
	& \le 
	2 \left\{ 
	\bigl| 
	\mbox{\boldmath $ v$} _{\varepsilon ,\lambda }(s) 
	\bigr|_{\mbox{\scriptsize \boldmath $ H$}_0}^2
	+ \bigl| 
	\pi (m_0)
	\bigr|^2 \big( |\Omega |+|\Gamma | \bigr) 
	\right\} 
	\nonumber \\
	& \le 
	2 \left\{ M_4
	+ \bigl| 
	\pi (m_0)
	\bigr|^2 \big( |\Omega |+|\Gamma | \bigr) 
	\right\}=:\tilde{M}_5
	\label{pi2}
\end{align}
for a.a.\ $s \in (0,T)$. 
We integrate the resultant \eqref{sq} 
over $(0,t)$ with respect to $s$. 
From Lemma~3.2, there exists a positive constant $M_5$, 
depending only on $c_3$, $c_4$, $T$, $|\Omega |$, $|\Gamma |$, $M_1$, $M_3$, $M_4$,
$\tilde{M}_5$, and 
$|\mbox{\boldmath $ f$}|_{L^2(0,T;\mbox{\scriptsize \boldmath $ H$}_0)}$ 
and independent of $\varepsilon \in (0,1/4]$ and $\lambda \in (0,1]$, such that \eqref{m4b} 
holds. 
\hfill $\Box$ \\

Now, recalling \eqref{mu} and 
using the fact that $\mbox{\boldmath $ v$}_{\varepsilon ,\lambda }'(t)$, 
$\partial \varphi (\mbox{\boldmath $ v$}_{\varepsilon ,\lambda }(t)) \in \mbox{\boldmath $ H$}_0$ for a.a.\ 
$t \in (0,T)$, we obtain 
\begin{gather}
	m \bigl( \mbox{\boldmath $ \mu $}_{\varepsilon ,\lambda } (t) \bigr) 
	= m \bigl(  \mbox{\boldmath $ \beta $}_\lambda
	\bigl (\mbox{\boldmath $ u$}_{\varepsilon, \lambda } (t) \bigr) 
	+ \varepsilon \mbox{\boldmath $ \pi $}
	\bigl (\mbox{\boldmath $ u$}_{\varepsilon, \lambda } (t) \bigr) 
	-\mbox{\boldmath $ f$}(t) \bigr)
	\quad {\rm for~a.a.\ } t \in (0,T).
	\label{ome}
\end{gather}

\paragraph{Lemma 3.4.}
{\it There exist two positive constants $M_6$ and $M_7$,
independent of $\varepsilon \in (0,1/4]$ and $\lambda \in (0,1]$, such that}
\begin{gather*}
	\int_{0}^{t}\bigl| 
	m \bigl( \mbox{\boldmath $ \mu $}_{\varepsilon ,\lambda } (s)\bigr) 
	\bigr|^2 ds
	\le M_6, 
%	\label{m5a}
	\\
	\int_{0}^{t} \bigl | \mbox{\boldmath $ \mu $}_{\varepsilon,\lambda } (s)
	\bigr |_{\mbox{\scriptsize \boldmath $ V$} }^2ds \le M_7
	\nonumber 
%	\label{m5b}
\end{gather*}
{\it for all $t \in [0,T]$. }

\paragraph{Proof.} 
From \eqref{mean} and \eqref{ome},
we have 
\begin{align*}
	\bigl| 
	m \bigl( \mbox{\boldmath $ \mu $}_{\varepsilon ,\lambda } (s) \bigr) 
	\bigr|^2
	& \le 
	\frac{6}{ \bigl(  |\Omega |+|\Gamma |\bigr)^2 }
	\left\{ 
	\bigl| 
	\beta _\lambda 
	\bigl( u_{\varepsilon, \lambda } (s) \bigr ) \bigr|^2_{L^1(\Omega )} 
	+ \varepsilon ^2 
	\bigl| 
	\pi \bigl( u_{\varepsilon, \lambda}(s) \bigr ) 
	\bigr| ^2_{L^1(\Omega )}
	+
	|\Omega |
	\bigl| f(s) \bigr|^2_{H} 
	\right.
	\\
	 & \left. {}+
	 \bigl|  
	\beta _{\lambda} \bigl( u_{\Gamma, \varepsilon,\lambda } (s) \bigr )
	\bigr|^2_{L^1(\Gamma)} +  \bigl|  
	\pi \bigl( u_{\Gamma ,\varepsilon,\lambda } (s) \bigr )
	\bigr|^2_{L^1(\Gamma  )} 
	+
	|\Gamma | \bigl|  
	f_\Gamma  (s) 
	\bigr|^2_{H_\Gamma}
	\right\}
\end{align*} 
for a.a.\ $s \in (0,T)$. 
Then, by integrating over $(0,T)$, 
it follows that there is a positive constant $M_6$, 
depending only on $T$, $|\Omega |$, $|\Gamma |$, $M_4$, $M_5$,   
$|\pi (m_0)|$, and 
$|\mbox{\boldmath $ f$}|_{L^2(0,T;\mbox{\scriptsize \boldmath $ H$}_0)}$, such that 
the first estimate holds. 
Next, from the {P}oincare inequality (see, e.g., \cite[Lemma~A]{CF15}), we know that 
there exists a positive constant $c_p$ such that 
\begin{equation} 
	c_p |\mbox{\boldmath $ z$}|_{\mbox{\scriptsize \boldmath $ V$}}^2 \le 
	|\mbox{\boldmath $ z$}|_{\mbox{\scriptsize \boldmath $ V$}_0}^2
	\quad {\rm for~all~}\mbox{\boldmath $ z$} \in \mbox{\boldmath $ V$}_0.
	\label{poin}
\end{equation} 
Therefore, 
\begin{align*}
	\int_{0}^{t} 
	\bigl | \mbox{\boldmath $ \mu $}_{\varepsilon,\lambda  }(s)
	\bigr |_{\mbox{\scriptsize \boldmath $ V$} } ^2 ds
	& \le 2 \int_{0}^{t} 
	\bigl | \mbox{\boldmath $ P$} \mbox{\boldmath $ \mu $}_{\varepsilon,\lambda  }(s)
	\bigr |_{\mbox{\scriptsize \boldmath $ V$} }^2 ds 
	+ 2 \int_{0}^{t} 
	\bigl| m \bigl( \mbox{\boldmath $ \mu $}_{\varepsilon ,\lambda } (s) \bigr) \mbox{\boldmath $ 1$}
	\bigr |_{\mbox{\scriptsize \boldmath $ V$} } ^2ds
	\\ 
	& \le \frac{2}{c_p} \int_{0}^{t} 
	\bigl | \mbox{\boldmath $ P$} \mbox{\boldmath $ \mu $}_{\varepsilon,\lambda } (s)
	\bigr |_{\mbox{\scriptsize \boldmath $ V$}_0 }^2 ds
	+ 2 \bigl( |\Omega | +|\Gamma | \bigr) 
	\int_{0}^{t}\bigl| 
	m \bigl( \mbox{\boldmath $ \mu $}_{\varepsilon ,\lambda } (s)\bigr) 
	\bigr|^2 ds \\
	& \le \frac{2M_3}{c_p}+ 2 \bigl( |\Omega | +|\Gamma | \bigr) M_6=:M_7
\end{align*} 
for all $t \in [0,T]$. \hfill $\Box$

\paragraph{Lemma 3.5.}
{\it There exists a positive constant $M_8$,
independent of $\varepsilon \in (0,1/4]$ and $\lambda \in (0,1]$, such that}
\begin{equation}
	\int_{0}^{t} 
	\bigl| \mbox{\boldmath $ \beta $}_\lambda 
	\bigl (\mbox{\boldmath $ u$}_{\varepsilon,\lambda }(s) \bigr) 
	\bigr|_{\mbox{\scriptsize \boldmath $ H$}}^2 ds 
	 \le M_8 \label{m6}
\end{equation}
{\it for all $t\in [0,T]$.}

\paragraph{Proof.} We test \eqref{mu} by 
$\mbox{\boldmath $ \beta $}_\lambda  (\mbox{\boldmath $ u$}_{\varepsilon,\lambda }(s)) 
\in \mbox{\boldmath $ V$}$, then
\begin{align*}
	& \varepsilon \bigl( \partial \varphi \bigl( 
	\mbox{\boldmath $ v$}_{\varepsilon ,\lambda } (s) \bigr), 
	\mbox{\boldmath $ \beta $}_\lambda 
	\bigl( 
	\mbox{\boldmath $ u$}_{\varepsilon ,\lambda }(s) 
	\bigr) 
	\bigr) _{\mbox{\scriptsize \boldmath $ H$}}
	+
	\bigl| \mbox{\boldmath $ \beta $}_\lambda 
	\bigl( 
	\mbox{\boldmath $ u$}_{\varepsilon ,\lambda } 
	(s) \bigr)\bigr|_{\mbox{\scriptsize \boldmath $ H$}}^2 
	\\
	& \le  
	 \bigl( \mbox{\boldmath $ \mu $}_{\varepsilon ,\lambda } (s)
	- \lambda 
	\mbox{\boldmath $ v $}_{\varepsilon ,\lambda }'(s) 
	- \varepsilon  \mbox{\boldmath $ \pi $} \bigl( 
	\mbox{\boldmath $ u $}_{\varepsilon ,\lambda } (s) \bigr)
	+ \mbox{\boldmath $ f $}(s), 
	\mbox{\boldmath $ \beta $}_\lambda \bigl( 
	\mbox{\boldmath $ u$}_{\varepsilon ,\lambda }(s) \bigr) 
	\bigr) _{\mbox{\scriptsize \boldmath $ H$}}
\end{align*} 
for a.a.\ $s \in (0,T)$. 
Here, from the monotonicity of $\beta _\lambda $, we have 
\begin{align*}
	& \varepsilon \bigl( \partial \varphi \bigl( 
	\mbox{\boldmath $ v$}_{\varepsilon ,\lambda } (s) \bigr), 
	\mbox{\boldmath $ \beta $}_\lambda 
	\bigl( 
	\mbox{\boldmath $ u$}_{\varepsilon ,\lambda }(s) 
	\bigr) 
	\bigr) _{\mbox{\scriptsize \boldmath $ H$}} 
	\\
	& = \varepsilon \int_{\Omega }^{} 
	(-\Delta v_{\varepsilon,\lambda  }(s)) \beta _\lambda \bigl( 
	u_{\varepsilon ,\lambda }(s) \bigr) dx 
	+ \varepsilon \int_{\Gamma }^{} \partial _{{\mbox{\scriptsize \boldmath $ \nu $}} } v_{\varepsilon ,\lambda }(s) 
	\beta_\lambda  \bigl(u_{\Gamma,\varepsilon ,\lambda  }(s) \bigr)d \Gamma \\
	& \quad {} + 
	\varepsilon \int_{\Gamma }^{} (-\Delta_\Gamma v_{\Gamma ,\varepsilon ,\lambda } (s)) 
	\beta _\lambda \bigl( u_{\Gamma,\varepsilon ,\lambda } (s)\bigr) d\Gamma \\
	& = \varepsilon 
	\int_{\Omega }^{} \beta '_\lambda
	\bigl( u_{\varepsilon,\lambda } (s) \bigr)
	\bigl | \nabla v_{\varepsilon,\lambda  }(s) \bigr |^2 dx 
	+ \int_{\Gamma }^{} \beta '_\lambda
	\bigl( u_{\Gamma, \varepsilon,\lambda } (s) \bigr)
	\bigl | \nabla_\Gamma v_{\Gamma,\varepsilon,\lambda } (s) \bigr |^2 d\Gamma \\
	& \ge 0
\end{align*} 
for a.a.\ $s \in (0,T)$, where we have used the fact 
that $\beta _\lambda 
(u_{\varepsilon ,\lambda }(s))_{|_\Gamma }
=\beta _\lambda (u_{\Gamma,\varepsilon ,\lambda}(s))$ a.e.\ on $\Gamma$. 
Moreover, by the Young inequality, 
\begin{align*} 
	& \bigl( \mbox{\boldmath $ \mu $}_{\varepsilon ,\lambda } (s)
	- \lambda 
	\mbox{\boldmath $ v $}_{\varepsilon ,\lambda }'(s) 
	- \varepsilon  \mbox{\boldmath $ \pi $} \bigl( 
	\mbox{\boldmath $ u $}_{\varepsilon ,\lambda } (s) \bigr) 
	+  \mbox{\boldmath $ f $}(s), 
	\mbox{\boldmath $ \beta $}_\lambda 
	\bigl( 
	\mbox{\boldmath $ u$}_{\varepsilon ,\lambda }(s) 
	\bigr) 
	\bigr) _{\mbox{\scriptsize \boldmath $ H$}}
	\\
	& \le 
	\frac{1}{2} \bigl| \mbox{\boldmath $ \beta $}_\lambda 
	\bigl( 
	\mbox{\boldmath $ u$}_{\varepsilon ,\lambda } 
	(s) \bigr)\bigr|_{\mbox{\scriptsize \boldmath $ H$}}^2
	+ 2 \bigl| \mbox{\boldmath $ \mu $}_{\varepsilon ,\lambda } (s) 
	\bigr|_{\mbox{\scriptsize \boldmath $ H$}}^2
	+ 2 \lambda ^2 
	\bigl| \mbox{\boldmath $ v $}_{\varepsilon ,\lambda }'(s) 
	\bigr|_{\mbox{\scriptsize \boldmath $ H$}_0}^2 
	+ 2 \varepsilon ^2
	\bigl| \mbox{\boldmath $ \pi $} \bigl( 
	\mbox{\boldmath $ u $}_{\varepsilon ,\lambda } (s) \bigr)
	\bigr|_{\mbox{\scriptsize \boldmath $ H$}}^2
	+ 2
	\bigl| \mbox{\boldmath $ f $}(s) 
	\bigr|_{\mbox{\scriptsize \boldmath $ H$}_0}^2
\end{align*} 
for a.a.\ $s \in (0,T)$. 
Let us combineing the above inequalities, and integrateing in over $(0,T)$ 
with respect to $s$.
Then, recalling \eqref{prim} and \eqref{pi2}, and using Lemmas~3.2 and 3.4,  
we obtain
\begin{align*}
	\int_{0}^{t} \bigl| 
	\mbox{\boldmath $ \beta $}_\lambda \bigl( 
	\mbox{\boldmath $ u$}_{\varepsilon ,\lambda }(s) \bigr) \bigr|_{\mbox{\scriptsize \boldmath $H $}} ^2 ds 
	& \le  
	4 \int_{0}^{T} \bigl| \mbox{\boldmath $ \mu $}_{\varepsilon ,\lambda } (s) 
	\bigr|_{\mbox{\scriptsize \boldmath $ H$}}^2 ds 
	+ 4 \lambda ^2 \int_{0}^{T} 
	\bigl| \mbox{\boldmath $ v $}_{\varepsilon ,\lambda }'(s) 
	\bigr|_{\mbox{\scriptsize \boldmath $ H$}_0}^2 ds \\
	& \quad {}
	+ 4  \varepsilon ^2 \int_{0}^{T}
	\bigl| \mbox{\boldmath $ \pi $} \bigl( 
	\mbox{\boldmath $ u $}_{\varepsilon ,\lambda } (s) \bigr)
	\bigr|_{\mbox{\scriptsize \boldmath $ H$}}^2
	ds
	+ 4 \int_{0}^{T}
	\bigl| \mbox{\boldmath $ f $}(s) 
	\bigr|_{\mbox{\scriptsize \boldmath $ H$}_0}^2ds
	\\
	& \le 
	4 M_7
	+ 4 M_3  + 4 T \tilde{M}_5
	+ 4 
	\bigl| \mbox{\boldmath $ f $}
	\bigr|_{L^2(0,T;\mbox{\scriptsize \boldmath $ H$}_0)}^2=:M_8
\end{align*} 
for all $t \in [0,T]$. \hfill $\Box$

%%%%% Section 3.2. %%%%%
\subsection{Uniform estimates for approximate solutions of (P)$_\varepsilon$}

In this subsection, based on the result of Proposition~2.1, 
we obtain uniform estimates for the solutions of (P)$_\varepsilon$. 
Actually, by obtaining additional uniform estimates (see \cite{CF15}) and 
passing to the limit 
in the approximate problem (P)$_{\varepsilon ,\lambda }$ 
as $\lambda \to 0$, Proposition~2.1 can be proved. 
Thus, for each $\varepsilon \in (0,1/4]$,
there exists a triplet $(\mbox{\boldmath $ v$}_\varepsilon, 
\mbox{\boldmath $ \mu  $}_\varepsilon, \mbox{\boldmath $ \xi $}_\varepsilon )$ 
of
\begin{gather*}
	\mbox{\boldmath $ v$}_\varepsilon 
	\in H^1(0,T;\mbox{\boldmath $ V$}_0^*) \cap L^\infty (0,T;\mbox{\boldmath $ V$}_0) \cap L^2(0,T;\mbox{\boldmath $ W$}), \\
	\mbox{\boldmath $ \mu $}_\varepsilon 
	\in L^2(0,T;\mbox{\boldmath $ V$}), \quad 
	\mbox{\boldmath $ \xi $}_\varepsilon 
	\in L^2(0,T;\mbox{\boldmath $ H$}), \quad 
	\mbox{\boldmath $ \xi $}_\varepsilon 
	\in \mbox{\boldmath $ \beta $}(\mbox{\boldmath $ u$}_\varepsilon )
	\quad {\rm in~} L^2(0,T;\mbox{\boldmath $ H$})
\end{gather*}
that satisfy \eqref{weak1}, \eqref{weak2}, and 
$\mbox{\boldmath $ v$}_\varepsilon (0)=\mbox{\boldmath $ v$}_{0}$ 
in $\mbox{\boldmath $ H$}_0$. Here, we also have 
\begin{equation*}
	\mbox{\boldmath $ u$}_\varepsilon 
	\in H^1(0,T;\mbox{\boldmath $ V$}^*) \cap L^\infty (0,T;\mbox{\boldmath $ V$}) 
	\cap L^2(0,T;\mbox{\boldmath $ W$}),
\end{equation*}
with 
the relation 
$\mbox{\boldmath $ u$}_\varepsilon =\mbox{\boldmath $ v$}_\varepsilon +m_0\mbox{\boldmath $ 1$}$. 
On taking the limit 
$\lambda \to 0$, we obtain the same kind of uniform estimates as obtained 
in the previous lemmas, namely
\begin{gather}
	\int_{0}^{t} 
	\bigl| 
	\mbox{\boldmath $ v$}_{\varepsilon} '(s)
	| _{\mbox{\scriptsize \boldmath $ V$}_0^*}^2 ds 
	\le M_3, 
	\label{e1} \\
	\varepsilon 
	\bigl|\mbox{\boldmath $ v$}_{\varepsilon}(t) 
	\bigr|_{\mbox{\scriptsize \boldmath $ V$}_0}^2  
	\le M_3, 
	\label{e2}\\
%	\int_{0}^{t} \bigl| \mbox{\boldmath $ P$} 
%	 \mbox{\boldmath $ \mu $}_{\varepsilon}(s) 
%	 \bigr|_{\mbox{\scriptsize \boldmath $ V$}_0}^2 ds 
%	 \le M_3,
%	\label{e3} \\
	\bigl| 
	\mbox{\boldmath $ u $}_{\varepsilon}(t) 
	\bigr|_{\mbox{\scriptsize \boldmath $ H$}}^2
	\le M_4, 
	\label{e4} \\
	\quad \bigl| 
	\mbox{\boldmath $ v $}_{\varepsilon}(t) 
	\bigr|_{\mbox{\scriptsize \boldmath $ H$}_0}^2
	\le M_4,
	\label{e5}
	\\
	\int_{0}^{t} \bigl | \mbox{\boldmath $ \mu $}_{\varepsilon} (s)
	\bigr |_{\mbox{\scriptsize \boldmath $ V$} }^2ds \le M_7,
	\label{e6} \\
	\int_{0}^{t} 
	\bigl| 
	\mbox{\boldmath $ \xi $}_{\varepsilon}(s) 
	\bigr|_{\mbox{\scriptsize \boldmath $ H$}}^2 ds 
	 \le M_8 
	 \label{e7}
\end{gather}
for all $t \in [0,T]$. Using these estimates, we can prove the main theorem.

%%%%% Section 4 %%%%%
\section{Proof of the main theorem}

In this subsection, we conclude the proof of the main theorem 
by passing to the limit in the approximate problem (P)$_\varepsilon $ 
as $\varepsilon \to 0$. 

\paragraph{Proof of Theorem 2.1.}
By using the estimates \eqref{e1}--\eqref{e7} stated in the previous section, 
there 
exist a subsequence $\{ \varepsilon _k \}_{k \in \mathbb{N}}$ with 
$\varepsilon _k \to 0$ as $k \to +\infty $
and some 
limit functions $\mbox{\boldmath $ u$} \in H^1(0,T;\mbox{\boldmath $ V$}^*) \cap L^\infty (0,T;\mbox{\boldmath $ H$})$, 
$\mbox{\boldmath $ v$} \in H^1(0,T;\mbox{\boldmath $ V$}_0^*) \cap L^\infty (0,T;\mbox{\boldmath $ H$}_0)$, 
$\mbox{\boldmath $ \mu $}\in L^2(0,T;\mbox{\boldmath $ V$})$, and 
$\mbox{\boldmath $ \xi $} \in L^2(0,T;\mbox{\boldmath $ H$})$ such that 
\begin{gather} 
	\mbox{\boldmath $ v$}_{\varepsilon_k} 
	\to \mbox{\boldmath $ v$} \quad \mbox{weakly star in } 
	L^\infty (0,T;\mbox{\boldmath $ H$}_0), 
	\label{c1}
	\\
	\mbox{\boldmath $ u$}_{\varepsilon_k} 
	\to \mbox{\boldmath $ u$}=\mbox{\boldmath $ v$}+m_0 \mbox{\boldmath $ 1$} \quad \mbox{weakly star in } 
	L^\infty (0,T;\mbox{\boldmath $ H$}), 
	\label{c2}
	\\
	\varepsilon _k \mbox{\boldmath $ v$}_{\varepsilon_k} \to \mbox{\boldmath $ 0$}
	\quad \mbox{in } 
	L^\infty (0,T;\mbox{\boldmath $ V$}_0), 
	\label{c3}
	\\
	\mbox{\boldmath $ v$}_{\varepsilon_k}' \to \mbox{\boldmath $ v$}' 
	\quad \mbox{weakly in } 
	L^2(0,T;\mbox{\boldmath $ V$}_0^*),
	\label{c4}
	\\
	\mbox{\boldmath $ u$}_{\varepsilon_k}' \to \mbox{\boldmath $ u$}' 
	\quad \mbox{weakly in } 
	L^2(0,T;\mbox{\boldmath $ V$}^*),
	\label{c5}
	\\
	\mbox{\boldmath $ \mu $}_{\varepsilon_k} \to \mbox{\boldmath $ \mu $} 
	\quad \mbox{weakly in } 
	L^2(0,T;\mbox{\boldmath $ V$}), 
	\label{c6}
	\\
%	\mbox{\boldmath $ P$}\mbox{\boldmath $ \mu $}_{\varepsilon_k} 
%	\to \mbox{\boldmath $ P$}\mbox{\boldmath $ \mu $} 
%	\quad \mbox{weakly in } 
%	L^2(0,T;\mbox{\boldmath $ V$}_0), 
%	\label{c7}
%	\\
	\mbox{\boldmath $ \xi $}_{\varepsilon_k} 
	\to \mbox{\boldmath $ \xi $} \quad \mbox{weakly in } 
	L^2(0,T;\mbox{\boldmath $ H$}),
	\label{c8}
	\\
	\varepsilon_k \mbox{\boldmath $\pi $} (\mbox{\boldmath $ u$}_{\varepsilon_k}) 
	\to \mbox{\boldmath $ 0$} \quad \mbox{in } 
	L^\infty (0,T;\mbox{\boldmath $ H$})
	\label{c9}
\end{gather} 
as $k \to +\infty$. 
From \eqref{c1}, \eqref{c2}, \eqref{c4}, and \eqref{c5}, 
the well-known {A}scoli--{A}rzela theorem
(see, e.g., \cite[Section~8, Corollary~4]{Sim87}) 
gives
\begin{gather} 
	\mbox{\boldmath $ v$}_{\varepsilon_k} \to \mbox{\boldmath $ v$} \quad \mbox{in } 
	C\bigl( [0,T];\mbox{\boldmath $ V$}_0^* \bigr),
	\label{c10} \\
	\mbox{\boldmath $ u$}_{\varepsilon_k} \to \mbox{\boldmath $ u$} \quad \mbox{in } 
	C\bigl( [0,T];\mbox{\boldmath $ V$}^* \bigr)
	\label{c11}
\end{gather}
as $k \to +\infty$. 
Now, integrating \eqref{weak2} over $(0,T)$, we find
\begin{align} 
	\int_{0}^{T} \bigl( 
	\mbox{\boldmath $ \mu $}_{\varepsilon_k} (t),\mbox{\boldmath $ \eta $}(t)
	\bigr)_{\mbox{\scriptsize \boldmath $ H$}} dt
	& = \varepsilon_k \int_{0}^{T} 
	a \bigl(  \mbox{\boldmath $ v$}_{\varepsilon_k} (t), 
	\mbox{\boldmath $ \eta $}(t) \bigr) dt
	+ \int_{0}^{T} 
	\bigl( 
	 \mbox{\boldmath $ \xi $}_{\varepsilon_k} (t), 
	 \mbox{\boldmath $ \eta $}(t) 
	 \bigr) _{\mbox{\scriptsize \boldmath $ H$}}
	 dt 
	 \nonumber \\
	 & 
	+ \varepsilon_k \int_{0}^{T}
	\bigl( 
	\mbox{\boldmath $ \pi$} \bigl( 
	\mbox{\boldmath $ u$}_{\varepsilon_k} (t)\bigr), \mbox{\boldmath $ \eta $}(t)
	\bigr)_{\mbox{\scriptsize \boldmath $ H$}} 
	dt
	-
	\int_{0}^{T}
	\bigl( \mbox{\boldmath $ f$}(t),\mbox{\boldmath $ \eta $}(t)
	\bigr)_{\mbox{\scriptsize \boldmath $ H$}} dt
	\nonumber 
	\\
	& \quad \mbox{\rm for all }\mbox{\boldmath $ \eta $} \in L^2(0,T;\mbox{\boldmath $ V$}).
	\label{weakvari}
\end{align}
Letting $k \to \infty $, we can use \eqref{c6}, \eqref{c3}, \eqref{c8}, 
and \eqref{c9} to obtain 
\begin{equation*} 
	\int_{0}^{T} \bigl( 
	\mbox{\boldmath $ \mu $} (t),\mbox{\boldmath $ \eta $}(t)
	\bigr)_{\mbox{\scriptsize \boldmath $ H$}} dt
	= 
	\int_{0}^{T} 
	\bigl( 
	 \mbox{\boldmath $ \xi $} (t)-\mbox{\boldmath $ f$}(t),\mbox{\boldmath $ \eta $}(t)
	\bigr)_{\mbox{\scriptsize \boldmath $ H$}} dt
	\quad \mbox{\rm for all }\mbox{\boldmath $ \eta $} \in L^2(0,T;\mbox{\boldmath $ V$}),
\end{equation*}
namely, $\mbox{\boldmath $ \mu $}=\mbox{\boldmath $ \xi $}-\mbox{\boldmath $ f$}$ in 
$L^2(0,T;\mbox{\boldmath $ H$})$. Additionally, 
the information $\mbox{\boldmath $ \mu $}+\mbox{\boldmath $ f$} \in L^2(0,T;\mbox{\boldmath $ V$})$ gives 
the regularity of $\mbox{\boldmath $ \xi$} \in L^2(0,T;\mbox{\boldmath $ V$})$,
namely $\xi _\Gamma =\xi _{|_\Gamma }$ a.e.\ on $\Sigma $. 
Next, we take $\mbox{\boldmath $ \eta $}:=\mbox{\boldmath $ u$}_{\varepsilon _k} 
\in L^2(0,T;\mbox{\boldmath $ V$})$ in 
\eqref{weakvari}. Then, using the positivity
\begin{equation*} 
	\varepsilon_k  \int_{0}^{T} 
	a \bigl( \mbox{\boldmath $ v$}_{\varepsilon_k} (t), 
	\mbox{\boldmath $ u $}_{\varepsilon _k}(t) \bigr) dt 
	= \varepsilon_k  \int_{0}^{T} 
	\bigl| \mbox{\boldmath $ v$}_{\varepsilon_k} (t) 
	\bigr|_{\mbox{\scriptsize \boldmath $ V$}}^2 dt \ge 0, 
\end{equation*} 
and recalling \eqref{c2}, \eqref{c6}, \eqref{c9}, and \eqref{c11}, we have
\begin{align*} 
	\limsup _{k \to +\infty }
	\int_{0}^{T} 
	\bigl( 
	 \mbox{\boldmath $ \xi $}_{\varepsilon_k} (t), 
	 \mbox{\boldmath $ u $}_{\varepsilon _k}(t) 
	 \bigr) _{\mbox{\scriptsize \boldmath $ H$}}
	 dt 
	 & \le 
	\int_{0}^{T} \bigl \langle \mbox{\boldmath $ u $}(t),
	\mbox{\boldmath $ \mu $} (t)
	\bigr \rangle_{\mbox{\scriptsize \boldmath $ V$}^*, \mbox{\scriptsize \boldmath $ V$}} dt
	+ 
	\int_{0}^{T}
	\bigl( \mbox{\boldmath $ f$}(t),\mbox{\boldmath $ u $}(t)
	\bigr)_{\mbox{\scriptsize \boldmath $ H$}} dt \\
	& = \int_{0}^{T} \bigl( 
	 \mbox{\boldmath $ \xi $} (t), 
	 \mbox{\boldmath $ u $}(t) 
	 \bigr) _{\mbox{\scriptsize \boldmath $ H$}}
	 dt. 
\end{align*}
Thus, applying \cite[Proposition~2.2, p.~38]{Bar10} with
\eqref{c1} and \eqref{c8}, we can deduce that
\begin{gather*}
	\mbox{\boldmath $ \xi $} \in \mbox{\boldmath $ \beta $}(\mbox{\boldmath $ u$}) 
	\quad \mbox{a.e.\ in } Q
\end{gather*}
from the maximal monotonicity of $\beta $. Namely, 
we obtain 
$\xi \in \beta (u)$ a.e.\ in $\Omega $, 
$\xi _\Gamma \in \beta (u_\Gamma )$ a.e.\ on $\Gamma $.
Finally, integrating \eqref{weak1} over $(0,T)$ with respect to $t$ 
and letting $k \to + \infty $, we have
\begin{equation*} 
	\int_{0}^{T} 
	\bigl \langle \mbox{\boldmath $ v$}' (t),\mbox{\boldmath $ \eta $}(t) 
	\bigr \rangle _{\mbox{\scriptsize \boldmath $ V$}_0^*, 
	\mbox{\scriptsize \boldmath $ V$}_0}dt
	+ \int_{0}^{T} 
	a\bigl( \mbox{\boldmath $ P$}\mbox{\boldmath $ \mu $} (t),\mbox{\boldmath $ \eta $}(t) \bigr) 
	dt 
	=0
	\quad 
	\mbox{\rm for all } 
	\mbox{\boldmath $ \eta $} \in L^2(0,T;\mbox{\boldmath $ V$}_0).
\end{equation*} 
Next, using \eqref{dual}, we deduce 
\begin{equation*} 
	\mbox{\boldmath $ v$}'(t)+\mbox{\boldmath $ F$} \bigl( 
	\mbox{\boldmath $ P$} \mbox{\boldmath $ \mu $}(t) \bigr) = 
	\mbox{\boldmath $ 0$} 
	\quad {\rm in~} \mbox{\boldmath $ V$}_0^*, 
	\quad {\rm for~a.a.~}t \in (0,T),
\end{equation*}
namely 
\begin{equation} 
	\langle \mbox{\boldmath $ v$}' (t),\mbox{\boldmath $ z $}
	\bigr \rangle _{\mbox{\scriptsize \boldmath $ V$}_0^*, 
	\mbox{\scriptsize \boldmath $ V$}_0}
	+ 
	a\bigl( \mbox{\boldmath $ \xi $} (t),\mbox{\boldmath $ z $} \bigr) 
	= 
	a\bigl( \mbox{\boldmath $ f $} (t),\mbox{\boldmath $ z $} \bigr) 
	\quad {\rm for~all}~\mbox{\boldmath $ z$} \in \mbox{\boldmath $ V$}_0.
	\label{www}
\end{equation}
Therefore, recalling \eqref{ax} and \cite[Remark~2]{CF15}, we finally obtain
\begin{align*} 
	\bigl \langle u' (t),z 
	\bigr \rangle _{V^*,V}
	& + \langle u'_{\Gamma} (t),z_\Gamma 
	\bigr \rangle _{V_\Gamma ^*,V_\Gamma }
	+ \int_{\Omega }^{} \nabla \xi (t) \cdot \nabla z dx 
	+ \int_{\Gamma }^{} \nabla_\Gamma  \xi _{\Gamma}(t) 
	\cdot \nabla _\Gamma z_\Gamma 
	d \Gamma 
	\\
	& 
	= \int_{\Omega }^{} g(t) z dx 
	+ \int_{\Gamma }^{} g_\Gamma (t) z_\Gamma 
	d \Gamma 
	\quad {\rm for~all~}\mbox{\boldmath $ z$}=(z,z_\Gamma ) 
	\in \mbox{\boldmath $ V$},
\end{align*}
for a.a.\ $t \in (0,T)$, with $u(0)=u_0$ a.e.\ in $\Omega $ and $u_\Gamma (0)=u_{0\Gamma }$ a.e.\ on $\Gamma $. 
Thus, it turns out that the pair $(\mbox{\boldmath $ u$}, \mbox{\boldmath $ \xi$})$ 
is a weak solution of (P). 

%%%%% Section 5 %%%%%
\section{Continuous dependence}

In this section, we prove the continuous dependence of the data. 
This theorem also guarantees the uniqueness of the component $\mbox{\boldmath $ u$}$ in 
the solution.

\paragraph{Proof of Theorem 2.2.} For $i=1,2$, let 
$(\mbox{\boldmath $ u$}^{(i)}, \mbox{\boldmath $ \xi $}^{(i)})$
be a weak solution of {\rm (P)} corresponding to the data 
($\mbox{\boldmath $ f$}^{(i)}$, $\mbox{\boldmath $ u$}_{0}^{(i)}$). 
Set $m_0^*:=m(\mbox{\boldmath $ u$}_{0}^{(1)})=m(\mbox{\boldmath $ u$}_{0}^{(2)})$. 
From the weak formulation \eqref{www} of the {S}tefan problem (P) for 
$\mbox{\boldmath $ v$}^{(i)}=\mbox{\boldmath $ u$}^{(i)}
-m_0^*\mbox{\boldmath $ 1$}$, we obtain
\begin{gather} 
	\bigl \langle (\mbox{\boldmath $ v$}^{(i)})'(s),\mbox{\boldmath $ z$} 
	\bigr \rangle _{\mbox{\scriptsize \boldmath $ V$}_0^*, 
	\mbox{\scriptsize \boldmath $ V$}_0}
	+ a \bigl( \mbox{\boldmath $ \xi $}^{(i)}(s),\mbox{\boldmath $ z$} \bigr) 
	= a \bigl( \mbox{\boldmath $ f $}^{(i)}(s),\mbox{\boldmath $ z$} \bigr) 
	\quad \mbox{\rm for all } \mbox{\boldmath $ z$} \in \mbox{\boldmath $ V$}_0,
	\label{we}
\end{gather}
for a.a.\ $s\in (0,T)$. 
Now, $\mbox{\boldmath $ P$}\mbox{\boldmath $ \mbox{\boldmath $ \xi $}$}^{(i)}(s)
=\mbox{\boldmath $ \xi $}^{(i)}(s)-m(\mbox{\boldmath $ \xi $}^{i}(s))\mbox{\boldmath $ 1$} 
\in \mbox{\boldmath $ V$}_0$. 
Therefore, we see that 
\begin{align}
	a \bigl( \mbox{\boldmath $ \xi $}^{(i)}(s), 
	 \mbox{\boldmath $ F$}^{-1}
	\mbox{\boldmath $ z$}_0
	\bigr) 
	& = a \bigl( \mbox{\boldmath $ F$}^{-1}
	\mbox{\boldmath $ z$}_0, \mbox{\boldmath $ P$}\mbox{\boldmath $ \xi $}^{(i)}(s) \bigr) \nonumber \\
	& = \bigl \langle 
	\mbox{\boldmath $ z$}_0,  \mbox{\boldmath $ P$}\mbox{\boldmath $ \xi $}^{(i)}(s)
	\bigr \rangle _{\mbox{\scriptsize \boldmath $ V$}_0^*,\mbox{\scriptsize \boldmath $ V$}_0} \nonumber \\
	& = \bigl( \mbox{\boldmath $ z$}_0,  \mbox{\boldmath $ P$}\mbox{\boldmath $ \xi $}^{(i)}(s)
	\bigr)_{\mbox{\scriptsize \boldmath $ H$}_0} 
	\nonumber \\
	& = \bigl( \mbox{\boldmath $ z$}_0, \mbox{\boldmath $ \xi $}^{(i)}(s)
	\bigr)_{\mbox{\scriptsize \boldmath $ H$}} 
	\quad {\rm for~all}~ \mbox{\boldmath $ z$}_0 \in \mbox{\boldmath $ H$}_0 \subset \mbox{\boldmath $ V$}_0^*.
\end{align}
Thus, 
we take the difference between equation \eqref{we} when $i=1$ and when $i=2$.  
Then, we set
$\mbox{\boldmath $ z$}:=\mbox{\boldmath $ F$}^{-1}
(\mbox{\boldmath $ v$}^{(1)}(s)-\mbox{\boldmath $ v$}^{(2)}(s)) 
\in \mbox{\boldmath $ V$}_0$,
and use \eqref{dual} to obtain
\begin{align} 
	& \frac{1}{2} \frac{d}{ds} 
	\bigl| \mbox{\boldmath $ v$}^{(1)}(s)-  \mbox{\boldmath $ v$}^{(2)}(s) 
	\bigr|_{\mbox{\scriptsize \boldmath $ V$}_0^*}^2 
	+ 
	\bigl( 
	\mbox{\boldmath $ v$}^{(1)}(s)-  \mbox{\boldmath $ v$}^{(2)}(s),
	\mbox{\boldmath $ \xi $}^{(1)}(s)-  \mbox{\boldmath $ \xi $}^{(2)}(s)
	\bigr)_{\mbox{\scriptsize \boldmath $ H$}}
	\nonumber \\
	& \le \bigl \langle
	\mbox{\boldmath $ v$}^{(1)}(s) -  \mbox{\boldmath $ v$}^{(2)}(s),
	\mbox{\boldmath $ f$}^{(1)}(s) -\mbox{\boldmath $ f$}^{(2)}(s)  
	\bigr \rangle_{\mbox{\scriptsize \boldmath $ V$}^*_0, 
	\mbox{\scriptsize \boldmath $ V$}_0} \nonumber  \\
	& \le \frac{1}{2} 
	\bigl|
	\mbox{\boldmath $ v$}^{(1)}(s) -  \mbox{\boldmath $ v$}^{(2)}(s)
	\bigr|_{\mbox{\scriptsize \boldmath $ V$}_0^*} ^2
	+
	\frac{1}{2}
	\bigl| 
	\mbox{\boldmath $ f$}^{(1)}(s) -\mbox{\boldmath $ f$}^{(2)}(s)  
	\bigr |_{\mbox{\scriptsize \boldmath $ V$}_0} ^2
	\label{pr2}
\end{align} 
for a.a.\ $s \in (0,T)$. 
Now, in the variational formulation \eqref{ax} used to construct 
$\mbox{\boldmath $ f$}$, taking 
$\mbox{\boldmath $ z$}:=\mbox{\boldmath $ f$}$, we obtain 
\begin{equation*} 
	\bigl| \mbox{\boldmath $ f$} 
	\bigr|_{\mbox{\scriptsize \boldmath $ V$}_0}^2
	\le 
	\frac{1}{2c_p} \bigl| \mbox{\boldmath $ g$} 
	\bigr|_{\mbox{\scriptsize \boldmath $ H$}_0}^2
	+\frac{c_p}{2}\bigl| \mbox{\boldmath $ f$} 
	\bigr|_{\mbox{\scriptsize \boldmath $ H$}_0}^2
\end{equation*} 
a.e.\ in $(0,T)$. 
Therefore, by using the {Y}oung and {P}oincare \eqref{poin} inequalities, 
\begin{align} 
	& \frac{1}{2} \bigl| \mbox{\boldmath $ f$}^{(1)}(s) -\mbox{\boldmath $ f$}^{(2)}(s)
	\bigr|_{\mbox{\scriptsize \boldmath $ V$}_0}^2 \nonumber \\
	& \le \frac{1}{2} \bigl| \mbox{\boldmath $ f$}^{(1)}(s) -\mbox{\boldmath $ f$}^{(2)}(s)
	\bigr|_{\mbox{\scriptsize \boldmath $ V$}}^2 
	\nonumber \\
	& \le  
	\frac{1}{2} \bigl| \mbox{\boldmath $ f$}^{(1)}(s) -\mbox{\boldmath $ f$}^{(2)}(s)
	\bigr|_{\mbox{\scriptsize \boldmath $ V$}}^2 
	+
	\frac{1}{2}
	\left\{  
	\bigl| \mbox{\boldmath $ f$}^{(1)}(s) -\mbox{\boldmath $ f$}^{(2)}(s)
	\bigr|_{\mbox{\scriptsize \boldmath $ V$}}^2 
	-\bigl| \mbox{\boldmath $ f$}^{(1)}(s) -\mbox{\boldmath $ f$}^{(2)}(s)
	\bigr|_{\mbox{\scriptsize \boldmath $ H$}_0}^2
	\right\} 
	\nonumber \\
	& \le \frac{1}{c_p}
	\bigl| \mbox{\boldmath $ f$}^{(1)}(s) -\mbox{\boldmath $ f$}^{(2)}(s)
	\bigr|_{\mbox{\scriptsize \boldmath $ V$}_0}^2 
	-
	\frac{1}{2} 
	\bigl| \mbox{\boldmath $ f$}^{(1)}(s) -\mbox{\boldmath $ f$}^{(2)}(s)
	\bigr|_{\mbox{\scriptsize \boldmath $ H$}_0}^2
	\nonumber \\
	& \le 
	\frac{1}{2c_p^2} 
	\bigl| \mbox{\boldmath $ g$}^{(1)}(s) -\mbox{\boldmath $ g$}^{(2)}(s)
	\bigr|_{\mbox{\scriptsize \boldmath $ H$}_0}^2
	\label{contikey}
\end{align} 
for a.a.\ $s \in (0,T)$. 
Then, from the monotonicity of $\beta $, 
we can apply the {G}ronwall inequality 
to deduce \eqref{conti1} for some positive constant $C$ that depends only on $c_p$ and 
$T$. Here, we must take care to ensure that 
$\mbox{\boldmath $ v$}^{(1)}-\mbox{\boldmath $ v$}^{(2)}=\mbox{\boldmath $ u$}^{(1)}-\mbox{\boldmath $ u$}^{(2)}$. 
If $\beta $ is a Lipschitz continuous function with 
{L}ipschitz constant $c_\beta>0$, then 
$\mbox{\boldmath $ \xi $}=(\beta (u),\beta (u_\Gamma ))$ and, with the help of the 
monotonicity of $\beta $, 
we have 
\begin{equation} 
	\bigl| \mbox{\boldmath $ \xi $}^{(1)}(s) 
	-  \mbox{\boldmath $ \xi $}^{(2)}(s) 
	\bigr|_{\mbox{\scriptsize \boldmath $ H$}}^2
	\le c_\beta \bigl( \mbox{\boldmath $ u $}^{(1)}(s) 
	-  \mbox{\boldmath $ u $}^{(2)}(s), 
	\mbox{\boldmath $ \xi $}^{(1)}(s) 
	-  \mbox{\boldmath $ \xi $}^{(2)}(s) 
	\bigr)_{\mbox{\scriptsize \boldmath $ H$}}.
	\label{pr3}
\end{equation} 
Thus, combining \eqref{pr2}, \eqref{contikey}, \eqref{pr3}, and \eqref{conti1}, and integrating 
the resultant quantity over $(0,t)$ with respect to $s$, we obtain
\begin{align*} 
	& 
	\bigl| \mbox{\boldmath $ u$}^{(1)}(t)-  \mbox{\boldmath $ u$}^{(2)}(t) 
	\bigr|_{\mbox{\scriptsize \boldmath $ V$}_0^*}^2 
	+ 
	\frac{2}{c_\beta }
	\int_{0}^{t}
	\bigl| 
	\mbox{\boldmath $ \xi $}^{(1)}(s)-  \mbox{\boldmath $ \xi $}^{(2)}(s)
	\bigr|_{\mbox{\scriptsize \boldmath $ H$}}^2ds
	\nonumber \\
	& \le 
	\bigl| \mbox{\boldmath $ u$}^{(1)}_0 -  \mbox{\boldmath $ v$}^{(2)}_0 
	\bigr|_{\mbox{\scriptsize \boldmath $ V$}_0^*}^2 
	+
	\int_{0}^{t}
	\bigl|
	\mbox{\boldmath $ u$}^{(1)}(s) -  \mbox{\boldmath $ u$}^{(2)}(s)
	\bigr|_{\mbox{\scriptsize \boldmath $ V$}_0^*} ^2ds
	+
	\frac{1}{c_p^2}
	\int_{0}^{t}
	\bigl| 
	\mbox{\boldmath $ g$}^{(1)}(s) -\mbox{\boldmath $ g$}^{(2)}(s)  
	\bigr |_{\mbox{\scriptsize \boldmath $ H$}_0} ^2
	ds
	\\
	&  
	\le  C ( 1+T )  \left\{ 
	\bigl| \mbox{\boldmath $ u$}^{(1)}_0-  \mbox{\boldmath $ u$}^{(2)}_0 
	\bigr|_{\mbox{\scriptsize \boldmath $ V$}_0^*}^2 
	+
	\int_{0}^{T} 
	\bigl| \mbox{\boldmath $ g$}^{(1)}(s ) -\mbox{\boldmath $ g$}^{(2)}(s ) 
	\bigr|_{\mbox{\scriptsize \boldmath $ H$}_0}^2 ds
	\right\}
\end{align*}
for all $t \in [0,T]$. That is, we have obtained \eqref{conti2}. \hfill $\Box$

\section*{Acknowledgments}

The author wishes to express his heartfelt gratitude to 
Professors {G}oro {A}kagi and {U}lisse {S}tefanelli, 
who kindly gave him the opportunity of exchange visits
supported by the JSPS--CNR bilateral joint research project on
\emph{{I}nnovative {V}ariational {M}ethods for {E}volution {E}quations}. 
The author is also indebted to Professor 
Pierluigi Colli, who kindly gave him the opportunity for fruitful discussions 
that aided in obtaining the results of this paper. 
The author was supported by JSPS KAKENHI 
Grant-in-Aid for Scientific Research(C), Grant Number 26400164.

\end{document}